\documentclass[12pt]{article}

\usepackage{latexsym}
\usepackage{amsmath}
\usepackage{amsbsy}
\usepackage{amsthm}
\usepackage{amssymb}
\usepackage{amsfonts}
\usepackage{dsfont}
\usepackage{mathrsfs}

\newtheorem{proposition}{Proposition}[section]
\newtheorem{theorem}[proposition]{Theorem}
\newtheorem{corollary}[proposition]{Corollary}
\newtheorem{lemma}[proposition]{Lemma}

\numberwithin{equation}{section}

\newcommand{\nc}{\newcommand}
\nc{\R}{{\mathds R}}
\nc{\N}{{\mathds N}}
\nc{\Z}{{\mathds Z}}
\nc{\BP}{\mathds{P}}
\nc{\BE}{\mathds{E}}
\nc{\BQ}{\mathds{Q}}
\nc{\BX}{\mathds{X}}
\nc{\I}{\mathds{1}}
\nc{\cA}{{\mathscr A}}
\nc{\cN}{{\mathscr N}}
\nc{\cX}{{\mathscr X}}
\nc{\cY}{{\mathscr Y}}
\nc{\cF}{{\mathscr F}}
\nc{\cB}{{\mathscr B}}
\nc{\cD}{{\mathscr D}}
\nc{\cH}{{\mathscr H}}
\nc{\cS}{{\mathscr S}}
\nc{\bH}{{\mathbf H}}
\nc{\bN}{{\mathbf N}}
\nc{\dint}{{\rm d}}

\DeclareMathOperator{\dom}{{dom}}
\DeclareMathOperator{\Var}{{\mathrm Var}}

\begin{document}

\author{G\"unter Last\thanks{Karlsruhe Institute of Technology,
    Institut f\"ur Stochastik, Kaiserstra{\ss}e 89, D-76128 Karlsruhe,
    Germany. Email: \texttt{guenter.last@kit.edu}} 
}

\title{Stochastic analysis for Poisson processes}
\date{\today}

\maketitle

\begin{abstract}
\noindent
This survey is a preliminary version of a chapter
of the forthcoming book \cite{PeccReitz15}.
The paper develops some basic theory for
the stochastic analysis of Poisson process on a general 
$\sigma$-finite measure space. After giving some
fundamental definitions and properties
(as the multivariate Mecke equation) the paper 
presents the Fock space representation of 
square-integrable functions of a Poisson process 
in terms of iterated difference operators. 
This is followed by the introduction of multivariate stochastic
Wiener-It\^o integrals and the discussion of their basic properties.
The paper then proceeds with proving the chaos expansion
of square-integrable Poisson functionals, and defining and discussing
Malliavin operators. Further topics are products of
Wiener-It\^o integrals and Mehler's
formula for the inverse of the Ornstein-Uhlenbeck generator
based on a dynamic thinning procedure. The survey
concludes with covariance identities, the
Poincar\'e inequality and the FKG-inequality.
\end{abstract}


\smallskip


\section{Basic properties of a Poisson process}

Let $(\BX, \mathscr X)$ be a measurable space. The idea of a
{\em point process} with {\em state space} $\BX$ is that of a random 
countable subset of $\BX$,
defined over a fixed probability space $(\Omega,\cA,\mathds P)$.
It is both convenient and mathematically fruitful to define
a point process as a random element $\eta$ in the space
$\bN_\sigma(\BX)\equiv \bN_\sigma$ of all $\sigma$-finite measure $\chi$
on $\BX$ such that $\chi(B)\in \Z_+\cup\{\infty\}$ for all $B\in \mathscr X$.
To do so, we equip $\bN_\sigma$ with the smallest $\sigma$-field 
$\mathscr{N}_\sigma(\BX)\equiv \mathscr{N}_\sigma$ of subsets
of   $\bN_\sigma$ such that $\chi\mapsto \chi(B)$ is measurable for all $B\in \mathscr X$.
Then $\eta:\Omega\to \bN_\sigma$ is a point process if and only if
$\{\eta(B)=k\}\equiv\{\omega\in\Omega:\eta(\omega,B)=k\}\in \cA$
for all $B\in \cX$ and all $k\in \Z_+$. Here we write
$\eta(\omega,B)$ instead of the more clumsy $\eta(\omega)(B)$.
We wish to stress that the results of this survey
do not require special (topological) assumptions on the
state space.

The {\em Dirac measure} $\delta_x$ at the point $x\in \BX$ is the measure 
on $\BX$ defined by $\delta_x(B)=\mathds 1_B(x)$, where 
$\mathds 1_B$ is the indicator function of $B\in\mathscr X$.
If $X$ is a random element of $\BX$, then $\delta_X$ is a point
process on $\BX$. Suppose, more generally, that $X_1,\dots,X_m$
are independent random elements in $\BX$ with
distribution $\mathds Q$. Then
\begin{align}\label{e1.1.2}
\eta:=\delta_{X_1}+\dots+\delta_{X_m}
\end{align}
is a point process on $\BX$. Because
\begin{align*}
\mathds P(\eta(B)=k)=\binom{m}{k}\mathds Q(B)^k (1-\mathds Q(B))^{m-k},\quad k=0,\dots,m,
\end{align*}
$\eta$ is referred to a {\em binomial process} with
{\em sample size} $m$ and {\em sampling distribution} $\mathds Q$.
Taking an infinite sequence $X_1,X_2,\dots$ of independent random variables
with distribution $\mathds Q$ and
replacing in $\eqref{e1.1.2}$ the deterministic sample size $m$ by
an independent $\Z_+$-valued 
random variable $\kappa$ (and interpreting an empty sum
as null measure) yields a {\em mixed binomial process}.
Of particular interest is the case, where
$\kappa$ has a Poisson distribution with parameter $\lambda\ge 0$,
see also \eqref{e1.1.4} below.
It  is then easy to check that
\begin{align}\label{Laplace}
\mathds E\exp\bigg[-\int u(x)\eta(\dint x)\bigg]
=\exp\bigg[-\int (1-e^{-u(x)})\mu(\dint x)\bigg],
\end{align}
for any measurable function $u:\BX\to[0,\infty)$,
where $\mu:=\lambda \mathds Q$. It is convenient to write this as
\begin{align}\label{Laplace2}
\mathds E\exp[-\eta(u)]=\exp\big[-\mu(1-e^{-u})\big],
\end{align}
where $\nu(u)$ denotes the integral of a measurable function $u$
with respect to a measure $\nu$.
Clearly,
\begin{align}\label{intensitym}
\mu(B)=\mathds E\eta(B),\quad B\in\cX,
\end{align}
so that $\mu$ is the {\em intensity measure} of $\eta$.
The identity \eqref{Laplace2} or elementary probabilistic arguments show
that $\eta$ has {\em independent increments}, that is
the random variables $\eta(B_1),\dots,\eta(B_m)$ are stochastically independent
whenever $B_1,\ldots,B_m\in \mathscr X$ are pairwise disjoint.
Moreover, $\eta(B)$ has a Poisson distribution with parameter
$\mu(B)$, that is
\begin{align}\label{e1.1.4}
\mathds P(\eta(B)=k)=\frac{\mu(B)^k}{k!}\exp[-\mu(B)],
\quad k\in\Z_+.
\end{align}

Let $\mu$ be a $\sigma$-finite measure on $\BX$.
A {\em Poisson process} with intensity measure $\mu$ is a point process
$\eta$ on $\BX$ with independent increments such that \eqref{e1.1.4} holds,
where an expression of the form $\infty e^{-\infty}$ is interpreted as $0$.
It is easy to see that these two requirements determine the
distribution $\mathds P_{\eta}:=\mathds P(\eta\in\cdot)$ of a 
Poisson process $\eta$. We have seen above that a Poisson process exists
for a finite measure $\mu$. In the general case
it can be constructed as a countable sum of independent Poisson
processes, see \cite{Mecke,Kallenberg,LastPenrose} for more detail.
Equation \eqref{Laplace2} remains valid. Another consequence of this
construction is a representation of the form
\begin{align}\label{e1.1.6}
\eta=\sum^{\eta(\BX)}_{n=1}\delta_{X_n},
\end{align}
where $X_1,X_2,\dots$ are random elements in $\BX$.
This is one of the reasons why it is sufficient 
to work with a general $\sigma$-finite measure space $(\BX,\mathscr X,\mu)$ 
and to define a Poisson process
as a random element in the space $\bN_\sigma$ of $\sigma$-finite
measures on $(\BX, \mathscr X)$.

Let $\eta$ be a Poisson process with intensity measure $\mu$.
A classical and extremely useful formula by Mecke \cite{Mecke}  says that
\begin{align}\label{Mecke}
\BE\int h(\eta,x)\eta(\dint x)=\BE\int h(\eta+\delta_{x},x)\mu(\dint x)
\end{align}
for all measurable $h:\bN_\sigma\times \BX \rightarrow[0,\infty]$.
One can use the mixed binomial representation to prove
this result for finite Poisson processes.
An equivalent formulation is
\begin{align}\label{Mecke2}
\BE\int h(\eta-\delta_x,x)\eta(\dint x)=\BE\int h(\eta,x)\mu(\dint x)
\end{align}
for all measurable $h:\bN_\sigma\times \BX\rightarrow[0,\infty]$.
Although $\eta-\delta_x$ is in general a signed measure, we can use
\eqref{e1.1.6} to see that
\begin{align*}
\int h(\eta-\delta_x,x)\eta(\dint x)=\sum_i h\bigg(\sum_{j\ne i}\delta_{X_j},X_i\bigg)
\end{align*}
is almost surely well defined.
Both \eqref{Mecke} and \eqref{Mecke2} characterize the distribution of
a Poisson process with given intensity measure $\mu$.
 
Equation \eqref{Mecke} admits a useful generalization involving
multiple integration. To formulate this version we
consider, for $m\in\N$, the $m$-th power 
$(\BX^m,\mathscr X^m)$ of $(\BX,\mathscr X)$.
Let $\eta$ be given by \eqref{e1.1.6}. 
We define another point process $\eta^{(m)}$ on $\BX^m$ by
\begin{align}\label{mum}
\eta^{(m)}(C)=\sideset{}{^{\ne}}\sum_{i_1,\ldots,i_m\le\eta(\BX)}\mathds 1_C(X_{i_1},\dots,X_{i_m}),
\quad C\in \mathscr X^m,
\end{align}
where the superscript $\ne$ indicates summation over $m$-tuples with
pairwise different entries. The multivariate version of
\eqref{Mecke} (see e.g.\ \cite{LastPenrose}) says that
\begin{align}\label{Meckem}\notag
\mathds E\int &h(\eta,x_1,\dots,x_m)\eta^{(m)}(\dint(x_1,\dots,x_m))\\ 
&=\mathds E\int h(\eta+\delta_{x_1}+\dots+\delta_{x_m},x_1,\dots,x_m)
\mu^m(\dint(x_1,\ldots,x_m)),
\end{align}
for all measurable $h:\bN_\sigma\times \BX^m\rightarrow[0,\infty]$. 
In particular the {\em factorial moment measures} of $\eta$ are given by
\begin{align}\label{factorialmoment}
\BE \eta^{(m)}=\mu^m,\quad m\in \N.
\end{align}
Of course \eqref{Meckem} remains true for a measurable 
$h:\bN_\sigma\times \BX^m\rightarrow\R$
provided that the right-hand side is finite when replacing $h$ with $|h|$.

\section{Fock space representation}

In the remainder of this paper
we consider a Poisson process $\eta$ on $\BX$ with $\sigma$-finite intensity measure
$\mu$ and distribution $\mathds P_\eta$. 

In this and later sections the following {\em difference operators}
will play a crucial role.
For any $f\in \mathbf{F}(\bN_\sigma)$ (the set of all measurable
functions from $\bN_\sigma$ to $\R$) and $x\in \BX$ the function
$D_xf\in \mathbf{F}(\bN_\sigma)$ is defined by
\begin{align}\label{e1.2.1}
D_xf(\chi):=f(\chi+\delta_x)-f(\chi),\quad \chi\in \bN_\sigma.
\end{align}
Iterating this definition, for $n\ge 2$ and $(x_1,\dots,x_n)\in \BX^n$
we define a function
$D^{n}_{x_1,\dots,x_n}f\in \mathbf{F}(\bN_\sigma)$
inductively by
\begin{align}\label{e1.2.4}
D^{n}_{x_1,\dots,x_{n}}f:=D^1_{x_{1}}D^{n-1}_{x_2,\ldots,x_{n}}f,
\end{align}
where $D^1:=D$ and $D^0f = f$.  Note that
\begin{align}\label{e1.2.6}
D^n_{x_1,\ldots,x_n}f(\chi)=\sum_{J \subset \{1,2,\ldots,n\}}(-1)^{n-|J|}
f \Big(\chi+\sum_{j\in J}\delta_{x_j}\Big),
\end{align}
where  $|J|$ denotes the number of elements of $J$. This
shows that $D^n_{x_1,\dots,x_n}f$ is symmetric in $x_1,\ldots,x_n$ and that
$(x_1,\dots,x_n,\chi)\mapsto D^n_{x_1,\dots,x_n}f(\chi)$ is measurable.  
We define symmetric and measurable functions $T_nf$ on $\BX^n$ by
\begin{align}\label{e1.2.2}
T_n f (x_1,\dots,x_n) :=\mathds E D^n_{x_1,\dots,x_n} f(\eta),
\end{align}
and we set $T_0f:=\mathds E f(\eta)$, whenever these expectations are defined.
By $\langle\cdot,\cdot\rangle_n$ we denote the
scalar product in $L^2(\mu^n)$ and by $\|\cdot\|_n$ the associated norm.
Let $L^2_s(\mu^n)$ denote the symmetric functions in $L^2(\mu^n)$.
Our aim is to prove that
the linear mapping $f  \mapsto (T_n(f))_{n\geq 0}$
is an isometry from $L^2(\mathds P_\eta)$ into the {\em Fock space} 
given by the direct sum of the spaces 
$L^2_s(\mu^n)$, $n\ge 0$, (with $L^2$ norms scaled by $n!^{-1/2}$)
and with $L^2_s(\mu^0)$ interpreted as $\R$. In Section \ref{sectionchaos}
we will see that this mapping is surjective. The result (and its proof) is from
\cite{LaPe11} and can be seen as a crucial first step in the stochastic
analysis on Poisson spaces.

\begin{theorem}\label{t1.2.1} Let $f,g\in L^2(\mathds P_\eta)$. Then 
\begin{align}\label{e1.2.3}
\mathds E f(\eta)g(\eta)=(\mathds E f(\eta))(\mathds E g(\eta))+\sum^\infty_{n=1}\frac{1}{n!}
\langle T_nf,T_ng \rangle_n,
\end{align}
where the series converges absolutely.
\end{theorem}

We will prepare the proof with some lemmas.
Let $\cX_{0}$ be the system of all measurable 
$B\in\cX$ having $\mu(B)<\infty$.
Let $\mathbf{F}_0$ be the space of all bounded and measurable
functions $v:\BX\to[0,\infty)$ vanishing outside some $B\in\cX_{0}$.
Let $\mathbf{G}$ denote the space of all
(bounded and measurable) functions $g:\bN_\sigma\to \R$
of the form
\begin{align}\label{e1.2.5}
g(\chi)=a_1e^{-\chi(v_1)}+\ldots+a_ne^{-\chi(v_n)},
\end{align}
where $n\in\N$, $a_1,\ldots,a_n\in\R$ and $v_1,\ldots,v_n\in \mathbf{F}_0$.

\begin{lemma}\label{l1.2.1}
Relation \eqref{e1.2.3} holds for $f,g\in\mathbf{G}$.
\end{lemma}
{\em Proof:} By linearity it suffices to consider
functions $f$ and $g$ of the form
\begin{align*}
f(\chi)=\exp[-\chi(v)],\quad g(\chi)=\exp[-\chi(w)]
\end{align*}
for $v,w\in \mathbf{F}_0$. Then we have for $n\ge 1$
that
$$
D^nf(\chi)=\exp[-\chi(v)](e^{-v}-1)^{\otimes n},
$$
where
$(e^{-v}-1)^{\otimes n}(x_1,\dots,x_n):=\prod_{i=1}^n (e^{-v(x_i)}-1)$.
From \eqref{Laplace2} we obtain that
\begin{align}\label{e1.2.8}
T_nf=\exp[-\mu(1-e^{-v})](e^{-v}-1)^{\otimes n}.
\end{align}
Since $v\in\mathbf{F}_0$ it follows that
$T_nf\in L^2_s(\mu^n)$, $n\ge 0$.
Using \eqref{Laplace2} again, we obtain that
\begin{align}\label{e1.2.12}
\BE f(\eta)g(\eta)=\exp[-\mu(1-e^{-(v+w)})].
\end{align}
On the other hand we have from \eqref{e1.2.8} 
(putting $\mu^0(1):=1$) that
\begin{align*}
\sum_{n=0}^\infty&\frac{1}{n!}\langle T_nf,T_ng\rangle_n\\
&=\exp[-\mu(1-e^{-v})]\exp[-\mu(1-e^{-w})]
\sum_{n=0}^\infty\frac{1}{n!}\mu^n(((e^{-v}-1)(e^{-w}-1))^{\otimes n})\\
&=\exp[-\mu(2-e^{-v}-e^{-w})]
\exp[\mu((e^{-v}-1)(e^{-w}-1))].
\end{align*}
This equals the right-hand side of \eqref{e1.2.12}. \qed

\bigskip

To extend  \eqref{e1.2.3} to general 
$f,g\in L^2(\BP_\eta)$ we need two further lemmas.

\begin{lemma}\label{l1.2.3}\rm 
The set $\mathbf{G}$ is dense in $L^2(\BP_\eta)$.
\end{lemma}
{\sc Proof:} Let $\mathbf{W}$ be the space of all
bounded measurable $g:\bN_\sigma\to\R$ that
can be approximated in $L^2(\BP_\eta)$ by functions
in $\mathbf{G}$. This space is closed under monotone and
uniformly bounded convergence and contains
the constant functions. 
The space $\mathbf{G}$ is stable under multiplication
and we denote by $\cN'$ 
the smallest $\sigma$-field on $\bN_\sigma$ such that $\chi\mapsto h(\chi)$
is measurable for all $h\in\mathbf{G}$.
A functional version of the monotone class theorem 
(see e.g.\ Theorem I.21 in \cite{DeMeyer78}) implies that
$\mathbf{W}$ contains any bounded $\cN'$-measurable $g$.
On the other hand we have that
$$
\chi(C)=\lim_{t\to 0+}t^{-1}(1 - e^{-t\chi(C)}),\quad \chi\in \bN_\sigma,
$$
for any  $C\in \cX$.
Hence $\chi\mapsto\chi(C)$ is $\cN'$-measurable 
whenever $C\in\cX_0$. Since $\mu$ is $\sigma$-finite, 
for any $C\in \cX$ there is a monotone sequence
$C_k\in\cX_0$, $k\in\N$, with union $C$, so that
$\chi\mapsto\chi(C)$ is $\mathcal{N}'$-measurable.
Hence $\cN'=\cN_\sigma$ and it follows that $\mathbf{W}$ contains all
bounded measurable functions. But then $\mathbf{W}$
is clearly dense in $L^2(\BP_\eta)$ and the proof of
the lemma is complete.\qed

\begin{lemma}\label{lemsubs} Suppose that
$f,f^1,f^2,\ldots\in L^2(\BP_\eta)$ satisfy
$f^k\to f$ in $L^2(\BP_\eta)$ as $k\to\infty$, and that
$h: \bN_\sigma \to [0,1]$ is measurable. 
Let $n\in\N$, let $C\in\cX_0$ and set $B:=C^n$. 
Then
\begin{align}\label{013}
\lim_{k\to\infty}\int_B
\BE [| D^n_{x_1,\dots,x_n} f(\eta) - D^n_{x_1,\ldots,x_n} f^k(\eta)|
h(\eta)] \mu^n(\dint(x_1,\dots,x_n))=0.
\end{align}
\end{lemma}
{\sc Proof:} 
By \eqref{e1.2.6}, the relation \eqref{013} is
implied by the convergence
\begin{align}\label{014}
\lim_{k\to\infty}\int_B
\BE \Big[ \Big|f\Big(\eta+ \sum_{i=1}^m \delta_{x_i} \Big)-
f^k\Big(\eta+ \sum_{i=1}^m \delta_{x_i} \Big) \Big|
h(\eta) \Big]\mu^n(\dint(x_1,\ldots,x_n))=0
\end{align}
for all $m\in\{0,\ldots,n\}$. For $m=0$ this is obvious.
Assume $m\in\{1,\ldots,n\}$. Then the integral in \eqref{014}
equals
\begin{align*}
\mu(C)^{n-m}&\BE\bigg[\int_{C^m} 
\Big|f\Big(\eta+\sum_{i=1}^m \delta_{x_i} \Big)-
f^k\Big(\eta+ \sum_{i=1}^m \delta_{x_i} \Big) \Big|h(\eta)
\mu^m(\dint(x_1,\ldots,x_m))\bigg]\\
&=\mu(C)^{n-m}\BE \bigg[\int_{C^m}  |f(\eta)-f^k(\eta)|
h \Big(\eta - \sum_{i=1}^n \delta_ {x_i} \Big)
\eta^{(m)}(\dint(x_1, \ldots , x_m))\bigg]\\
& \le\mu(C)^{n-m}\BE[|f(\eta)-f^k(\eta)|\eta^{(m)}(C^m)],
\end{align*}
where we have used \eqref{Meckem} to get the equality.
By the Cauchy-Schwarz inequality the last expression is bounded
above by 
\begin{align*}
\mu(C)^{n-m}(\BE[(f(\eta)-f^k(\eta))^2])^{1/2}
(\BE[(\eta^{(m)}(C^m))^2])^{1/2}.
\end{align*}
Since the Poisson distribution has moments of all orders,
we obtain \eqref{014} and hence the lemma.\qed

\bigskip

{\em Proof of Theorem \ref{t1.2.1}:} By linearity and
the polarization identity
$$
4\langle u,v\rangle_n=
\langle u+v,u+v\rangle_n-\langle u-v,u-v\rangle_n
$$
it suffices to prove \eqref{e1.2.3} for $f=g \in L^2(\BP_\eta)$.
By Lemma \ref{l1.2.3} there are $f^k\in \mathbf{G}$, $k\in\N$, satisfying
$f^k\to f$ in $L^2(\BP_\eta)$ as $k\to\infty$.
By Lemma \ref{l1.2.1}, $Tf^k$, $k\in\N$, is a Cauchy sequence
in $\bH:=\R\oplus\oplus^\infty_{n=1}L^2_s(\mu^n)$. 
The direct sum of the scalar products $(n!)^{-1}\langle \cdot,\cdot\rangle_n$
makes $\bH$ a Hilbert space. 
Let $\tilde f=(\tilde f_n)\in \bH$ be the limit, that is
\begin{align}\label{0lim}
\lim_{k\to\infty}
\sum_{n=0}^\infty\frac{1}{n!}\|T_nf^k-\tilde f_n\|^2_n=0.
\end{align}
Taking the limit in the identity
$\BE[f^k(\eta)^2]=\langle Tf^k,Tf^k\rangle_\bH$
yields $\BE[f(\eta)^2]=\langle \tilde f,\tilde f\rangle_\bH$.
Equation \eqref{0lim} implies that $\tilde f_0=\BE[f(\eta)]=T_0f$.
It remains to show that for any $n\ge 1$,
\begin{align}\label{identify}
\tilde f_n=T_nf,\quad \mu^n\text{-a.e.}
\end{align}
Let $C\in\cX_0$ and $B:=C^n$. 
Let $\mu^n_B$ denote the restriction of the measure $\mu^n$ to $B$. 
By \eqref{0lim} 
$T_nf^k$ converges in $L^2(\mu^n_B)$ (and hence in $L^1(\mu^n_B)$) to
$\tilde f_n$, while by the definition \eqref{e1.2.2} of $T_n$, and the case 
$h \equiv 1$ of \eqref{014}, 
$T_nf^k$ converges in $L^1(\mu^n_B)$ to $T_nf$. Hence 
these $L^1(\BP)$ limits must be the same almost everywhere, so that 
$\tilde f_n=T_nf$  $\mu^n$-a.e.\ on $B$. Since $\mu$
is assumed $\sigma$-finite, this implies \eqref{identify}
and hence the theorem. \qed

\section{Multiple Wiener-It\^o integrals}

For $n\ge 1$ and $g\in L^1(\mu^n)$ we define
\begin{align}\label{e1.3.1}
I_n(g):=\sum_{J \subset [n]}
(-1)^{n-|J|}\iint g(x_1,\dots,x_n)\eta^{(|J|)} (\dint x_J)
\mu^{n-|J|}(\dint x_{J^c}),
\end{align}
where $[n]:=\{1,\dots,n\}$, $J^c:=[n]\setminus J$ and $x_J:=(x_j)_{j\in J}$.
If $J=\emptyset$, then the inner integral on the right-hand side
has to be interpreted as $\mu^n(g)$. (This is to say that $\eta^{(0)}(1):=1$.)
The multivariate Mecke equation \eqref{Meckem} implies that
all integrals in \eqref{e1.3.1} are finite and that
$\BE I_n(g)=0$.

Given functions $g_i: \BX \to \R$ for $ i = 1,\dots,n$,
the {\em tensor product} 
$\otimes_{i=1}^n g_i$  is the function from $\BX^n$ to $\R$
which maps each $(x_1,\dots,x_n)$ to $\prod_{i=1}^n g_i(x_i)$. 
When the functions $g_1,\dots,g_n$ are all the same
function $h$, we write $h^{\otimes n}$ for this tensor
product function. In this case the definition \eqref{e1.3.1}
simplifies to
\begin{align}\label{e1.3.5}
I_n(h^{\otimes n})=\sum_{k =0}^n\binom{n}{k}(-1)^{n-k}
\eta^{(k)}(h^{\otimes k})(\mu(h))^{n-k}.
\end{align}

Let $\Sigma_n$ denote the set of all permutations of 
$[n]$, and for $g \in \BX^n\to\R$ define the
{\em symmetrization} $\tilde g$ of $g$ by 
\begin{align}\label{e1.3.3}
\tilde g(x_1,\ldots,x_n):=\frac{1}{n!}\sum_{\pi \in \Sigma_n} 
g(x_{\pi(1)},\ldots,x_{\pi(n)}).
\end{align}

The following {\em isometry properties} of the operators
$I_n$ are crucial. The proof is
similar to the one of \cite[Theorem 3.1]{LaPeSchulThae14}
and is based on the product form \eqref{factorialmoment}
of the factorial moment measures and some combinatorial arguments.
For more information on the intimate relationships
between moments of Poisson integrals and the combinatorial
properties of partitions we refer to \cite{Surg84,TaPecc11,LaPeSchulThae14}.

\begin{lemma}\label{l1.3.1}
Let $g\in L^2(\mu^m)$ and $h\in L^2(\mu^n)$ for $m,n\ge 1$ and assume that
$\{g\ne 0\}\subset B^m$ and $\{h\ne 0\}\subset B^n$ for some
$B\in\cX_0$. Then
\begin{align}\label{e1.orth}
  \BE I_m(g)I_n(h)=\I\{m=n\}m!\langle \tilde{g}, \tilde{h} \rangle_m.
\end{align}
\end{lemma}
{\em Proof:}  We start with a combinatorial identity. Let $n\in\N$.
A {\em subpartition} of $[n]$ is a (possibly empty) family
$\sigma$ of non-empty pairwise disjoint subsets of $[n]$.
The cardinality of $\cup_{J\in\sigma}J$ is denoted by $\|\sigma\|$.
For $u\in \mathbf{F}(\BX^{n})$ we define
$u_\sigma:\BX^{|\sigma|+n-\|\sigma\|}\to\R$ by
identifying the arguments belonging to the same $J\in\sigma$.
(The arguments $x_1,\dots,x_{|\sigma|+n-\|\sigma\|}$ have to be inserted in the
order of occurrence.) Now we take $r,s\in\Z_+$
such that $r+s\ge 1$ and define  $\Sigma_{r,s}$ as the set of all
partititons of $\{1,\dots,r+s\}$ such that
$|J\cap \{1,\dots,r\}|\le 1$ and $|J\cap \{r+1,\dots,r+s\}|\le 1$
for all $J\in\sigma$. Let $u\in\mathbf{F}(\BX^{r+s})$.
It is easy to see that
\begin{align}\label{e1.3.2}\notag
\iint u(x_1,\dots,x_{r+s}) \eta^{(r)}(\dint(x_1,\dots,x_{r}))
&\eta^{(s)}(\dint (x_{r+1},\dots,x_{r+s}))\\
&=\sum_{\sigma\in \Sigma_{r,s}} \int u_\sigma\, \dint\eta^{(|\sigma|)},
\end{align}
provided that $\eta(\{u\ne 0\})<\infty$. (In the case $r=0$ the inner integral on
the left-hand side is interpreted as $1$.)

We next note that $g\in L^1(\mu^m)$ and $h\in L^1(\mu^n)$
and abbreviate $f:=g \otimes h$. Let $k:=m+n$,
$J_1:=[m]$ and $J_2:=\{m+1,\dots,m+n\}$.
The definition \eqref{e1.3.1} and Fubini's theorem imply that
\begin{equation}\label{e1.3.4}
\begin{split}
I_{m}(g)I_{n}(h)=
\sum_{I\subset[k]}&(-1)^{n-|I|}\iiint f(x_1,\dots,x_k)\\
&\eta^{(|I\cap J_1|)}(\dint x_{I\cap J_1})
\eta^{(|I\cap J_2|)}(\dint x_{I\cap J_2})
\mu^{n-|I|}(\dint x_{I^c}),
\end{split}
\end{equation}
where $I^c:=[k]\setminus I$ and $x_J:=(x_j)_{j\in J}$ for any $J\subset[k]$.
We now take the expectation of
\eqref{e1.3.4} and use Fubini's theorem (justified by our integrability
assumptions on $g$ and $h$). Thanks 
to \eqref{e1.3.2} and \eqref{factorialmoment}
we can compute the expectation of the inner two integrals to obtain that
\begin{align}\label{e1.3.6}
\BE I_{m}(g)I_{n}(h)=
\sum_{\sigma\in \Sigma^*_{m,n}}(-1)^{k-\|\sigma\|}\int f_\sigma\,
\dint\mu^{k-\|\sigma\|+|\sigma|},
\end{align}
where $\Sigma^*_{m,n}$ is the set of all
subpartititons $\sigma$ of $[k]$ such that
$|J\cap J_1|\le 1$ and $|J\cap J_2|\le 1$ for all $J\in\sigma$.
Let $\Sigma^{*,2}_{m,n}\subset \Sigma^{*}_{m,n}$ be the set of all
subpartititons of $[k]$ such that $|J|=2$ for all $J\in\sigma$.
For any $\pi\in\Sigma^{*,2}_{m,n}$ we let $\Sigma^{*}_{m,n}(\pi)$ denote the
set of all $\sigma\in\Sigma^{*}_{m,n}$ satisfying $\pi\subset \sigma$.
Note that $\pi\in\Sigma^{*}_{m,n}(\pi)$ and that for any $\sigma\in\Sigma^{*}_{m,n}$
there is a unique $\pi\in \Sigma^{*,2}_{m,n}$ such that 
$\sigma\in\Sigma^{*}_{m,n}(\pi)$. In this case 
\begin{align*}
\int f_\sigma\dint\mu^{k-\|\sigma\|+|\sigma|}
=\int f_\pi\dint\mu^{k-\|\pi\|},
\end{align*}
so that \eqref{e1.3.6} implies
\begin{align}\label{e1.3.8}
\BE I_{m}(g)I_{n}(h)=
\sum_{\pi\in \Sigma^{*,2}_{m,n}} \int f_\pi\dint\mu^{k-\|\pi\|}
\sum_{\sigma\in \Sigma^*_{m,n}(\pi)}(-1)^{k-\|\sigma\|}.
\end{align}
The inner sum comes to zero, except in the case where
$\|\pi\|=k$. Hence \eqref{e1.3.8} vanishes unless $m=n$.
In the latter case we have
\begin{align*}
\BE I_{m}(g)I_{n}(h)=
\sum_{\pi\in \Sigma^{*,2}_{m,m}:|\pi|=m} \int f_\pi\,\dint\mu^{m}
=m!\langle \tilde{g}, \tilde{h} \rangle_m,
\end{align*}
as asserted.\qed

\bigskip

Any $g\in L^2(\mu^m)$ is the $L^2$-limit of
a sequence $g_k\in L^2(\mu^m)$ satisfying the assumptions
of Lemma \ref{l1.3.1}. For instance we may take
$g_k:=\I_{(B_k)^m}g$, where $\mu(B_k)<\infty$ and
$B_k\uparrow \BX$ as $k\to\infty$. Therefore the isometry
\eqref{e1.orth} allows to extend the linear operator $I_m$ 
in a unique way to $L^2(\mu^m)$. It follows from the isometry that
$I_m (g) = I_m(\tilde g)$ for all $g \in L^2(\mu^m)$.
Moreover, \eqref{e1.orth} remains true for arbitrary $g\in L^2(\mu^m)$ and
$h\in L^2(\mu^n)$. It is convenient to set
$I_0(c):=c$ for $c\in\R$. When $m\ge 1$,
the random variable $I_m(g)$ is the ($m$-th order)
{\em Wiener-It\^o integral}  of $g\in L^2(\mu^m)$ with respect to
the {\em compensated Poisson process} $\hat\eta:=\eta-\mu$.
The reference to $\hat\eta$ comes from the explicit definition \eqref{e1.3.1}.
We note that $\hat\eta(B)$ is only defined for $B\in\cX_0$. 
In fact, $\{\hat\eta(B):B\in\cX_0\}$
is an {\em independent random measure} in the sense of \cite{Ito56}.
The explicit definition \eqref{e1.3.1} was noted in \cite{Surg84}.

Let  $g\in L^2(\mu)$ and $f\in L^2(\mu^n)$ for some $n\in\N$.
Sometimes it is useful to write $I_n(f)I_1(g)$
as a sum of stochastic integrals. The following result
from \cite{Kab75} shows how this can be done. For any
$j\in[n]$ we define a function $f\otimes^0_j g:\BX^n\rightarrow \R$
by
\begin{align}\label{e1.2.81}
f\otimes^0_j g(x_1,\dots,x_n):=f(x_1,\dots,x_n)g(x_j)
\end{align}
and a function $f\otimes^1_j g:\BX^{n-1}\rightarrow \R$
by
\begin{align}\label{e1.2.82}
f\otimes^1_jg(x_1,\dots,x_{j-1},x_{j+1},\dots,x_n)
:=\int f\otimes^0_j g(x_1,\dots,x_n)\mu(\dint x_j).
\end{align}
By the Cauchy-Schwarz inequality
the latter integrals are finite $\mu^{n-1}$-a.e.

\begin{proposition}\label{pproduct} Let $n\in\N$,
$f\in L^2(\mu^n)$ and $g\in L^2(\mu)$. Assume that
$f\otimes^0_j g\in L^2(\mu^n)$ for all $j\in[n]$.
Then $f\otimes^1_j g\in L^2(\mu^{n-1})$ for all $j\in[n]$
and
\begin{align}\label{e1.2.7}
I_n(f)I_1(g)=I_{n+1}(f\otimes g)+\sum^n_{j=1}I_n(f\otimes^0_j g)
+\sum^n_{j=1}I_{n-1}(f\otimes^1_j g),\quad \BP\text{-a.s.}
\end{align}
\end{proposition}
{\em Proof:} The first assertion is a quick consequence of
the Cauchy-Schwarz inequality, see \eqref{e1.2.77} below.

To prove \eqref{e1.2.7} we first assume in addition that
$f$ and $g$ satsify the assumptions of Lemma \ref{l1.3.1}.
We can then use 
\eqref{e1.3.4} to obtain that
\begin{align}\label{e1.2.71}
I_n(f)I_1(g)=A_1-A_2,
\end{align}
where
\begin{align*}
A_1:=\sum_{J \subset [n]}
(-1)^{n-|J|}\iiint f(x_1,\dots,x_n)g(y)\eta (\dint y)
\eta^{(|J|)} (\dint x_J)\mu^{n-|J|}(\dint x_{[n]\setminus J})
\end{align*}
and
\begin{align*}
A_2:=\sum_{J \subset [n]}
(-1)^{n-|J|}\iiint f(x_1,\dots,x_n)g(y)
\eta^{(|J|)} (\dint x_J)\mu^{n-|J|}(\dint x_{[n]\setminus J})\mu(\dint y).
\end{align*}
From  the definition \eqref{mum} of the factorial measures
we see that $A_1$ equals
\begin{align*}
\sum_{J \subset [n]}&
(-1)^{n-|J|}\iint f(x_1,\dots,x_n)g(x_{n+1})
\eta^{(|J|+1)} (\dint x_{J\cup\{n+1\}})\mu^{n-|J|}(\dint x_{[n]\setminus J})\\
+&
\sum_{J \subset [n]}
(-1)^{n-|J|}\iint \sum_{j\in J}f(x_1,\dots,x_n)g(x_j)
\eta^{(|J|)} (\dint x_{J})\mu^{n-|J|}(\dint x_{[n]\setminus J}).
\end{align*}
The first sum can be rewritten as a sum over all 
$J \subset [n+1]$ with $n+1\in J$. Moreover, it is easy to check
that the sum without this restriction gives $I_{n+1}(f\otimes g)$.
This yields (after some rearranging)
\begin{align*}
A_1=&I_{n+1}(f\otimes g)+A_2\\
&+\sum^n_{j=1}
\sum_{\substack{J \subset [n]\\j\in J}
}
(-1)^{n-|J|}\iint f(x_1,\dots,x_n)g(x_j)
\eta^{(|J|)} (\dint x_{J})\mu^{n-|J|}(\dint x_{[n]\setminus J}).
\end{align*}
Therefore we obtain from \eqref{e1.2.71} that
\begin{align*}
I_n(f)&I_1(g)=I_{n+1}(f\otimes g)+\sum^n_{j=1}I_n(f\otimes^0_j g)\\
&-\sum^n_{j=1}
\sum_{\substack{J \subset [n]\\j\notin J}
}
(-1)^{n-|J|}\iint f(x_1,\dots,x_n)g(x_j)
\eta^{(|J|)} (\dint x_{J})\mu^{n-|J|}(\dint x_{[n]\setminus J})
\end{align*}
and \eqref{e1.2.7} follows. 

In the general case we define, for $k\in\N$, $f_k:=\I_{(B_k)^n}f$
and $g_k:=\I_{B_k}g$, where $\mu(B_k)<\infty$ and
$B_k\uparrow \BX$ as $k\to\infty$. Then we have not only
$f_k\to f$ in $L^2(\mu^n)$ and $g_k\to g$ in $L^2(\mu)$, but also
$f_k\otimes^0_jg_k\to f\otimes^0_j g$ in $L^2(\mu^n)$ for
any $j\in[n]$.
We have already shown that
\begin{align}\label{e1.2.74}
I_n(f_k)I_1(g_k)=I_{n+1}(f_k\otimes g_k)+\sum^n_{j=1}I_n(f_k\otimes^0_j g_k)
+\sum^n_{j=1}I_{n-1}(f_k\otimes^1_j g_k).
\end{align}
By the triangle and the Cauchy-Schwarz inequality the left-hand
side tends to $I_n(f)I_1(g)$ in $L^1(\BP)$ as $k\to\infty$.
We show that the right-hand side converges in $L^2(\BP)$.
Indeed, the isometry \eqref{e1.orth} and the Minkowski inequality yield that
\begin{align*}
\big[&\BE \big(I_{n+1}(f\otimes g)-I_{n+1}(f_k\otimes g_k))^2\big]^{1/2}
=\big[\BE \big(I_{n+1}(f\otimes g-f_k\otimes g_k))^2]^{1/2}\\
&=\sqrt{(n+1)!} \big[\mu^{n+1}\big((f\otimes g-f_k\otimes g_k)^2\big)\big]^{1/2}\\
&\le \sqrt{(n+1)!} \big[\mu^{n+1}(\big(f\otimes (g-g_k))^2\big)\big]^{1/2}
+\sqrt{(n+1)!} \big[\mu^{n+1}(\big((f-f_k)\otimes g_k)^2\big)\big]^{1/2}\\
&=\sqrt{(n+1)!} \big[\mu^{n}(f^2)\mu\big((g-g_k)^2\big)\big]^{1/2}
+\sqrt{(n+1)!} \big[\mu^{n}\big((f-f_k)^2\big)\mu\big((g_k)^2\big)\big]^{1/2}.
\end{align*}
As $k\to\infty$, this tends to $0$.
The other terms in \eqref{e1.2.74} can be treated in a similar way. For instance,
\begin{align*}
\big[&\BE \big(I_{n-1}(f\otimes^1_j g)-I_{n-1}(f_k\otimes^1_j g_k))^2\big]^{1/2}\\
&=\sqrt{(n+1)!}\big[\mu^{n-1}\big((f\otimes^1_j g-f_k\otimes^1_j g_k)^2\big)\big]^{1/2}\\
&\le \sqrt{(n+1)!}\big[\mu^{n-1}\big((f\otimes^1_j (g-g_k))^2\big)\big]^{1/2}
+\sqrt{(n+1)!}\big[\mu^{n-1}\big((f-f_k)\otimes^1_j g_k)^2\big)\big]^{1/2}.
\end{align*}
By the Cauchy-Schwarz inequality,
\begin{align}\label{e1.2.77}
&\mu^{n-1}\big((f\otimes^1_j (g-g_k))^2\big)\\ \notag
&\;\le\iint f(x_1,\dots,x_n)^2\mu(\dint x_j)\mu\big((g-g_k)^2\big)
\mu^{(n-1)}(\dint (x_1,\dots,x_{j-1},x_{j+1},\dots,x_n))\\ \notag
&\;=\mu^{n}\big(f^2)\mu\big((g-g_k)^2\big).
\end{align}
This tends to $0$ as $k\to\infty$. Similarly we get that
$\mu^{n-1}\big((f-f_k)\otimes^1_j g_k)^2\big)\to 0$.\qed

\bigskip

In Section \ref{secproduct} we will generalize Proposition \ref{pproduct}
to products $I_p(f)I_q(g)$,
where  $f\in L^2(\mu^p)$ and $g\in L^2(\mu^q)$ for $p,q\in\N$.

\section{The Wiener-It\^o chaos expansion}\label{sectionchaos}

A fundamental result of It\^o \cite{Ito56} and Wiener \cite{Wiener38} says
that every square integrable function of the Poisson process $\eta$ can be written
as an infinite series of orthogonal stochastic integrals.
Our aim is to prove the following explicit 
version of this {\em Wiener-It\^o chaos expansion}. Recall definition
\eqref{e1.2.2}.

\begin{theorem}\label{tchaos0} 
Let $f\in L^2(\BP_\eta)$. 
Then $T_nf\in L^2_s(\mu^n)$, $n\in\N$, and 
\begin{align}\label{chaos}
f(\eta)=\sum^\infty_{n=0}\frac{1}{n!}I_n(T_nf),
\end{align}
where the series converges in $L^2(\BP)$. 
Moreover, if $ g_n \in L^2_s(\mu^n)$ for
$n\in\Z_+$ satisfy 
$f(\eta)=\sum^\infty_{n=0}\frac{1}{n!}I_n(g_n)$ with
convergence in $L^2(\BP)$, then $g_0= \BE f(\eta)$ and
$g_n=T_nf$, $\mu^n$-a.e.\ on $\BX^n$, for all $n\in\N$.
\end{theorem}

For a homogeneous Poisson process on the real line,
the explicit chaos expansion \eqref{chaos} was proved in \cite{Ito88}.
The general case was formulated and proved in \cite{LaPe11}.
Stroock \cite{stroock87} has proved the counterpart of 
\eqref{chaos} for Brownian motion. Stroock's formula involves
iterated Malliavin derivatives 
and requires stronger integrability assumptions on $f(\eta)$.

Theorem \ref{tchaos0} and the isometry properties
\eqref{e1.orth} of stochastic integrals  show
that the isometry $f\mapsto (T_n(f))_{n\geq 0}$
is in fact a bijection from $L^2(\BP_\eta) $ onto the Fock space.
The following lemma is the key for the proof.

\begin{lemma}\label{l1.3.3} Let $f(\chi):=e^{-\chi(v)}$, $\chi\in\bN_\sigma(\BX)$,
where $v:\BX\to[0,\infty)$
is a measurable function vanishing outside a set $B\in\cX$ with $\mu(B)<\infty$.
Then \eqref{chaos} holds $\BP$-a.s.\ and in $L^2(\BP)$.
\end{lemma}
{\em Proof:} 
By \eqref{Laplace2} and \eqref{e1.2.8}
the right-hand side of \eqref{chaos} equals the formal sum
\begin{align}\label{e1.3.12}
I:=\exp[-\mu(1-e^{-v})]
+\exp[-\mu(1-e^{-v})]\sum^\infty_{n=1}\frac{1}{n!}
I_n((e^{-v}-1)^{\otimes n}).
\end{align}
Using the pathwise definition \eqref{e1.3.1} we obtain that almost surely
\begin{align}
I&=\exp[-\mu(1-e^{-v})]\sum^\infty_{n=0}\frac{1}{n!}
\sum^n_{k=0}\binom{n}{k}
\eta^{(k)}((e^{-v}-1)^{\otimes k})(\mu(1-e^{-v}))^{n-k}
\nonumber \\
&=\exp[-\mu(1-e^{-v})]
\sum^\infty_{k=0}\frac{1}{k!}\eta^{(k)}((e^{-v}-1)^{\otimes k})
\sum^\infty_{n=k}\frac{1}{(n-k)!}(\mu(1-e^{-v}))^{n-k}
\nonumber \\
&=\sum^{N}_{k=0}\frac{1}{k!}\eta^{(k)}((e^{-v}-1)^{\otimes k}),
\label{0825a}
\end{align}
where $N:=\eta(B)$. Writing
  $\delta_{X_1}+\dots+\delta_{X_N}$
for the restriction of $\eta$ to $B$, we have almost surely that
\begin{align*}
I= \sum_{J\subset  \{1,\dots,N\} }\prod_{i\in J}(e^{-v(X_i)}-1)
&=\prod^N_{i=1}e^{-v(X_i)}=e^{-\eta(v)},
\end{align*}
and hence \eqref{chaos} holds with almost sure convergence of the series.
To demonstrate that convergence also holds in $L^2(\BP)$,
let the partial sum $I(m)$ be given by the right hand side \eqref{e1.3.12} with  
the series terminated at $n=m$. Then since $\mu(1-e^{-v})$
is nonnegative and $|1-e^{-v(y)}| \le 1$ for all $y$,
a similar argument to \eqref{0825a} yields
\begin{align*}
|I(m)| & \leq  
\sum^{\min(N,m)}_{k=0}\frac{1}{k!} | \eta^{(k)}((e^{-v}-1)^{\otimes k}) | \\
 & \leq  
\sum^{N}_{k=0}\frac{N(N-1) \cdots (N-k+1)}{k!} 
= 2^N.
\end{align*}
Since $2^N$ has finite moments of all orders, by dominated
convergence the series
\eqref{e1.3.12} (and hence \eqref{chaos}) converges in $L^2(\BP)$. \qed

\bigskip

{\em Proof of Theorem \ref{tchaos0}:} Let $f \in L^2(\BP_\eta)$
and define $T_n f$ for $n\in\Z_+$ by \eqref{e1.2.2}. 
By \eqref{e1.orth} and  Theorem \ref{t1.2.1}, 
$$
\sum^\infty_{n=0}\BE\Big(\frac{1}{n!}I_n(T_nf)\Big)^2
=\sum^\infty_{n=0}\frac{1}{n!}\|T_nf\|_n^2 = \BE f(\eta)^2 <\infty.
$$
Hence the infinite series of orthogonal terms
$$
S:=\sum^\infty_{n=0}\frac{1}{n!}I_n(T_nf)
$$
converges in $L^2(\BP)$. Let $h \in\mathbf{G}$,  
where $\mathbf{G}$ was defined at \eqref{e1.2.5}. 
By Lemma \ref{l1.3.3} and linearity of $I_n(\cdot)$
the sum $\sum_{n=0}^\infty \frac{1}{n!} I_n(T_n h)$
converges in $L^2(\BP)$ to $ h(\eta)$.  
Using \eqref{e1.orth} followed by Theorem \ref{t1.2.1} yields
\begin{align*}
\BE (h(\eta) - S)^2=\sum^\infty_{n=0}\frac{1}{n!}\|T_nh - T_n f\|_{n}
= \BE(f(\eta) -h(\eta))^2.
\end{align*}
Hence if $\BE(f(\eta) - h(\eta))^2$ is small, then so is
$\BE(f(\eta) - S)^2$.
Since $\mathbf{G}$ dense in $L^2(\BP_\eta)$
by  Lemma \ref{l1.2.3}, it follows that  $f(\eta)=S$ almost
surely.

To prove the uniqueness, suppose that 
also $ g_n \in L^2_s(\mu^n)$ for
$n \in \Z_+$ are such that
$\sum^\infty_{n=0}\frac{1}{n!}I_n(g_n)$ 
converges in $L^2(\BP)$ to $f(\eta)$.
By taking expectations  we must have $g_0=\BE f(\eta)= T_0 f$.
For $n\ge 1$ and $h \in L^2_s(\mu^n)$, by \eqref{e1.orth} and \eqref{chaos}
we have
$$
\BE f(\eta) I_n(h) = \BE I_n(T_n f) I_n(h)= n!\langle T_n f, h \rangle_n 
$$ 
and similarly with $T_n f$ replaced by $g_n$, so that
$\langle T_nf -g_n, h \rangle_n =0$. Putting $h= T_nf-g_n$
gives $\| T_nf -g_n \|_n =0$ for each $n$, completing   
the proof of the theorem. \qed

\section{Malliavin operators}\label{s1.5}

For any $p\ge 0$ we denote by $L^p_\eta$ 
the space of all random variables $F\in L^p(\BP)$
such that $F=f(\eta)$ $\BP$-almost surely, for some 
$f\in\mathbf{F}(\bN_\sigma)$. Note that the space $L^p_\eta$ is a subset
of $L^p(\BP)$ while $L^p(\BP_\eta)$ is the space
of all measurable functions  $f\in\mathbf{F}(\bN_\sigma)$
satisfying $\int |f|^p\,\dint\BP_\eta=\BE |f(\eta)|^p<\infty$.
The {\em representative} $f$ of $F\in L^p(\BP)$ is
is $\BP_\eta$-a.e.\ uniquely defined  element of $L^p(\BP_\eta)$.
For $x\in \BX$ we can then define
the random variable $D_xF:=D_xf(\eta)$. More generally, we
define $D^n_{x_1,\dots,x_n}F:=D^n_{x_1,\dots,x_n}f(\eta)$ 
for any $n\in\N$ and $x_1,\dots,x_n\in \BX$.
The mapping $(\omega,x_1,\dots,x_n)\mapsto D^n_{x_1,\dots,x_n}F(\omega)$
is denoted by $D^nF$ (or by $DF$ in the case $n=1$).
The multivariate Mecke equation \eqref{Meckem} easily implies that
these definitions are $\BP\otimes\mu$-a.e.\ independent of
the choice of the representative.

By \eqref{chaos} any $F\in L^2_\eta$ can be written as 
\begin{align}\label{chaos2}
F=\BE F + \sum_{n=1}^\infty I_n(f_n),
\end{align}
where $f_n :=\frac{1}{n!}\BE D^nF$. In particular we obtain from
\eqref{e1.orth} (or directly from Theorem \ref{t1.2.1}) that
\begin{align}\label{e1.4.1}
\BE F^2=(\BE F)^2+\sum_{n=1}^\infty n! \|f_n\|_n^2.
\end{align}
We denote by $\dom D$ the set of all $F\in L^2_\eta$
satisfying 
\begin{align}\label{e1.4.2}
\sum_{n=1}^\infty n n! \|f_n\|_n^2 <\infty.
\end{align}
The following result is taken from \cite{LaPe11} and
generalizes Theorem 6.5 in \cite{Ito88}
(see also Theorem 6.2 in \cite{NuViv90}).
It shows that under the assumption \eqref{e1.4.2} the pathwise
defined difference operator $DF$ coincides with the
{\em Malliavin derivative} of $F$.
The space $\dom D$ is the {\em domain} of this operator.

\begin{theorem}\label{t1.4.1} Let $F\in L^2_\eta$ be given by
\eqref{chaos2}. Then $DF\in L^2(\BP\otimes\mu)$ iff $F\in\dom D$.
In this case we have $\BP$-a.s.\ and for $\mu$-a.e.\ $x\in \BX$ that
\begin{align}\label{Dx}
D_x F=\sum_{n=1}^\infty n I_{n-1}(f_n(x,\cdot)).
\end{align}
\end{theorem}

The proof Theorem \ref{t1.4.1} requires some preparations.
Since 
\begin{align*}
\int\Big(\sum^\infty_{n=1}n n!\|f_n(x,\cdot)\|^2_{n-1}\Big)\mu(\dint x)
=\sum^\infty_{n=1}n n!\int \|f_n\|^2_n,
\end{align*}
\eqref{e1.orth} implies that the infinite series
\begin{align}\label{e1.4.6}
D'_xF:=\sum^\infty_{n=1} nI_{n-1}f_n(x,\cdot)
\end{align}
converges in $L^2(\BP)$ for $\mu$-a.e.\ $x\in \BX$
provided that $F\in\dom D$.
By construction of the stochastic integrals we can assume that
$(\omega,x)\mapsto (I_{n-1}f_n(x,\cdot))(\omega)$ is measurable
for all $n\ge 1$. Therefore we can also assume that
the mapping $D'F$ given by $(\omega,x)\mapsto D'_xF(\omega)$
is measurable. We have just seen that
\begin{align}\label{e1.4.8}
\BE \int (D'_xF)^2\mu(\dint x)=\sum^\infty_{n=1}n n!\int \|f_n\|^2_n,\quad F\in\dom D.
\end{align}

Next we introduce an operator acting on random functions
that will turn out to be the {\em adjoint} of the difference operator $D$,
see Theorem \ref{t1.4.3}.
For $p\ge 0$ let $L^p_\eta(\BP\otimes \mu)$ denote  the set of all
$H\in L^p(\BP\otimes \mu)$ satisfying
$H(\omega,x)=h(\eta(\omega),x)$ for $\BP\otimes \mu$-a.e.\ $(\omega,x)$
for some {\em representative} $h\in \mathbf{F}(\bN_\sigma\times \BX)$.
For such a $H$ we have for $\mu$-a.e.\ $x$
that $H(x):=H(\cdot,x)\in L^2(\BP)$ and (by Theorem \ref{tchaos0})
\begin{align}\label{e1.4.32}
H(x)=\sum_{n=0}^\infty I_n(h_n(x,\cdot)), \quad \BP\text{-a.s.},
\end{align}
where $h_0(x):=\BE H(x)$ and
$h_n(x,x_1,\dots,x_n):=\frac{1}{n!}\BE D^n_{x_1,\dots,x_n}H(x)$.
We can then define the {\em Kabanov-Skorohod integral} 
\cite{Hitsuda72,Kab75,Skor75,KabSk75} of $H$, denoted
$\delta(H)$, by
\begin{align}\label{e1.4.11}
\delta(H):=\sum^\infty_{n=0} I_{n+1}(h_n),
\end{align}
which converges in $L^2(\BP)$ provided that
\begin{align}\label{domainS}
\sum^\infty_{n=0}(n+1)!\int \tilde h_n^2\dint\mu^{n+1}<\infty.
\end{align}
Here \begin{align}\label{e1.4.10}
\tilde h_n(x_1,\dots,x_{n+1}):=
\frac{1}{(n+1)!}\sum^{n+1}_{i=1}
\BE D^{n}_{x_1,\ldots,x_{i-1},x_{i+1},\dots,x_{n+1}}H(x_i)
\end{align}
is the symmetrization of $h_n$.
The set of all $H\in L^2_\eta(\BP\otimes\mu)$ satisfying
the latter assumption is the domain $\dom \delta$ of the operator $\delta$.

We continue with a preliminary version of Theorem \ref{t1.4.3}.

\begin{proposition}\label{p1.4.2} Let $F\in \dom D$. 
Let $H\in L^2_\eta(\BP\otimes \mu)$ be given by \eqref{e1.4.32}
and assume that
\begin{align}\label{e1.4.9}
\sum^\infty_{n=0}(n+1)!\int h_n^2d\mu^{n+1}<\infty.
\end{align}
Then
\begin{align}\label{adj2}
\BE\int (D'_xF) H(x)\mu(\dint x)=\BE F\delta(H).
\end{align}
\end{proposition}
{\em Proof:} 
Minkowski inequality implies \eqref{domainS} and hence $H\in\dom\delta$. 
Using  \eqref{e1.4.6} and \eqref{e1.4.32}
together with \eqref{e1.orth}, we obtain that
\begin{align*}
  \BE\int (D'_xF) H(x)\mu(\dint x)
=\int \bigg(\sum^\infty_{n=1} n! \langle f_n(x,\cdot),h_{n-1}(x,\cdot)\rangle_{n-1}
\bigg)\mu(\dint x),
\end{align*}
where the use of Fubini's theorem is justified by \eqref{e1.4.8},
the assumption on $H$ and the Cauchy-Schwarz inequality.
Swapping the order of summation and integration (to be justified soon)
we see that the last integral equals
\begin{align*}
\sum^\infty_{n=1}n! \langle f_n,h_{n-1}\rangle_n=
\sum^\infty_{n=1}n! \langle f_n,\tilde h_{n-1}\rangle_n,
\end{align*}
where we have used the fact that $f_n$ is a symmetric function.
By definition \eqref{e1.4.11} and \eqref{e1.orth}, 
the last series coincides with $\BE F\delta(H)$.
The above change of order is permitted since
\begin{align*}
  \sum^\infty_{n=1} n!\int |\langle f_n(x,\cdot),&h_{n-1}(x,\cdot)\rangle_{n-1}|\mu(\dint x)\\
&\le \sum^\infty_{n=1} n!\int \|f_n(x,\cdot)\|_{n-1}\|h_{n-1}(x,\cdot)\|_{n-1}\mu(\dint x)
\end{align*}
and the latter series is finite in view of the 
Cauchy-Schwarz inequality, the finiteness of \eqref{chaos2} 
and assumption \eqref{e1.4.9}. \qed

\bigskip

{\em Proof of Theorem \ref{t1.4.1}:} We need to show that
\begin{align}\label{e1.4.19}
DF=D'F, \quad \BP\otimes\mu\text{-a.e.}
\end{align}
First consider
the case with $f(\chi) = e^{-\chi(v)} $ with a measurable
$v:\BX\to[0,\infty)$ vanishing outside a set with finite $\mu$-measure.
Then $n!f_n = T_nf$ is given by \eqref{e1.2.8}.
Given $n \in \N$,
\begin{align*}
n \cdot n!
\int f_n^2 d \mu^n  = \frac{1}{(n-1)!} \exp[2 \mu(e^{-v}-1)] 
(\mu((e^{-v}-1)^2))^{n}
\end{align*}
which is summable in $n$, so \eqref{e1.4.2} holds in this case.
Also, in this case, $D_xf(\eta) = (e^{v(x)}-1) f(\eta)$ by
\eqref{e1.2.1}, while
$f_{n}(\cdot,x)=(e^{-v(x)} -1) n^{-1} f_{n-1}$
so that by \eqref{e1.4.6},
$$
D'_x f( \eta)=\sum_{n=1}^\infty(e^{-v(x)} -1) I_{n-1}(f_{n-1} ) = 
(e^{-v(x)} -1) f(\eta) 
$$
where the last inequality is from Lemma  \ref{l1.3.3} again.
Thus \eqref{e1.4.19} holds for $f$ of this form. 
By linearity this extends to all elements of $\mathbf{G}$. 

Let us now consider the general case. 
Choose $g_k\in \mathbf{G}$, $k\in\N$, such that $G_k:=g_k(\eta)\to F$ in $L^2(\BP)$
as $k\to\infty$, see Lemma \ref{l1.2.3}.
Let $H\in L^2_\eta(\BP_\eta\otimes\mu)$ have the representative
$h(\chi,x):=h'(\chi)\I_B(x)$, where $h'$ is as in
Lemma \ref{l1.3.3} and $B\in\cX_0$.  From
Lemma \ref{l1.3.3} it is easy to see that \eqref{e1.4.9}
holds. Therefore we obtain from Proposition \ref{p1.4.2} 
and the linearity of the operator $D'$ that
\begin{align}\label{adj7}
\BE\int (D'_xF -D'_xG_k)H(x)\mu(\dint x)=\BE (F-G_k)\delta(H)
\to 0\quad \text{as $k\to\infty$}.
\end{align}
On the other hand,
$$
\BE\int (D_xF-D_xG_k)H(x)\mu(\dint x)
= \int_B\BE [(D_xf(\eta)-D_xg_k(\eta))h'(\eta)]\mu(\dint x),
$$
and  by the case $n=1$ of Lemma \ref{lemsubs},
this tends to zero as $k \to \infty$.
Since $D'_xg_k=D_xg_k$ a.s.\ for $\mu$-a.e.\ $x$
we obtain from \eqref{adj7} that
\begin{align}\label{49}
\BE\int (D'_xf)h(\eta,x)\mu(\dint x)=\BE\int (D_xf(\eta))h(\eta,x)\mu(\dint x).
\end{align}
By Lemma \ref{l1.2.3}, the linear combinations
of the functions $h$ considered above are dense in
$L^2(\BP_\eta\otimes\mu)$, and by linearity 
\eqref{49} carries through to $h$ in this dense class of functions too, 
so we may conclude that the assertion \eqref{e1.4.19} holds.

It follows from \eqref{e1.4.8} and \eqref{e1.4.19} that $F\in\dom D$ 
implies $DF\in L^2_\eta(\BP\otimes\mu)$. The other implication was noticed
in \cite[Lemma 3.1]{PeTh13}. To prove it, we assume 
$DF\in L^2_\eta(\BP\otimes\mu)$ and apply the Fock space representation
\eqref{e1.2.3} to $\BE (D_xF)^2$ for $\mu$-a.e.\ $x$. This gives
\begin{align*}
\int \BE (D_xF)^2\mu(\dint x)
&=\sum^\infty_{n=0}\frac{1}{n!} \iint (\BE D^{n+1}_{x_1,\dots,x_n,x})^2
\mu^n(\dint(x_1,\dots,x_n))\mu(\dint x)\\
&=\sum^\infty_{n=0} (n+1)(n+1)!\|f_{n+1}\|_{n+1}^2
\end{align*}
and hence $F\in\dom D$.
\qed

The following duality relation (also referred to as
partial integration) shows that the operator $\delta$
is the adjoint of the difference operator $D$.
It is a special case of
Proposition 4.2 in \cite{NuViv90} applying to general
Fock spaces. 

\begin{theorem}\label{t1.4.3} Let $F\in \dom D$ and $H\in\dom \delta$.
Then
\begin{align}\label{partialint}
  \BE\int (D_xF) H(x)\mu(\dint x)=\BE F\delta(H).
\end{align}
\end{theorem}
{\em Proof:} We fix $F\in \dom D$.
Theorem \ref{t1.4.1} and Proposition \ref{p1.4.2} imply
that \eqref{partialint} holds if $H\in L^2_\eta(\BP\otimes\mu)$ 
satisfies the stronger assumption \eqref{e1.4.9}.
For any $m\in\N$ we define
\begin{align}\label{e1.4.20}
H^{(m)}(x):=\sum^m_{n=0} I_n(h_n(x,\cdot)),\quad x\in \BX.
\end{align}
Since $H^{(m)}$ satisfies \eqref{e1.4.9} we obtain that
\begin{align}\label{e1.4.21}
  \BE\int (D_xF) H^{(m)}(x)\mu(\dint x)=\BE F\delta(H^{(m)}).
\end{align}
From \eqref{e1.orth} we have
\begin{align*}
\int \BE(H(x)-H^{(m)}(x))^2\mu(\dint x)
&=\int\bigg(\sum^\infty_{n=m+1}n! \|h_n(x,\cdot)\|^2_n\bigg)\mu(\dint x)\\
&=\sum^\infty_{n=m+1}n! \|h_n\|^2_{n+1}.
\end{align*}
As $m\to\infty$ this tends to zero, since
\begin{align*}
\BE\int H(x)^2\mu(\dint x)=\int \BE(H(x))^2\mu(\dint x)=\sum^\infty_{n=0}n! \|h_n\|^2_{n+1}
\end{align*}
is finite. It follows that the left-hand side of
\eqref{e1.4.21} tends to the left-hand side of \eqref{partialint}.

To treat the right-hand side of \eqref{e1.4.21} we note that
\begin{align}\label{e1.4.24}
\BE \delta(H-H^{(m)})^2=\sum^\infty_{n=m+1}\BE (I_{n+1}(h_n))^2
=\sum^\infty_{n=m+1}(n+1)! \|\tilde{h}_n\|^2_{n+1}.
\end{align}
Since $H\in\dom\delta$ this tends to $0$ as $m\to\infty$. Therefore
$\BE(\delta(H)-\delta(H^{(m)}))^2\to 0$ and the
right-hand side of
\eqref{e1.4.21} tends to the right-hand side of \eqref{partialint}.\qed

\bigskip

We continue with a basic isometry property of the Kabanov-Skorohod integral.
In the present generality the result is in \cite{LaPeSchu14}.
A less general version is \cite[Proposition 6.5.4]{Priv09}.

\begin{theorem}\label{t1.4.7}
Let $H\in L^2_\eta(\BP\otimes \mu)$ be such that
\begin{align}\label{e1.4.30}
\BE \iint (D_yH(x))^2\mu(\dint x)\mu(\dint y)<\infty.
\end{align}
Then, $H\in\dom\delta$ and moreover
\begin{align}\label{e1.4.31}
\BE \delta(H)^2= \BE\int H(x)^2\mu(\dint x)
+\BE \iint D_yH(x)D_xH(y)\mu(\dint x)\mu(\dint y).
\end{align}
\end{theorem}
{\em Proof:} Suppose that $H$ is given as in \eqref{e1.4.32}.
Assumption \eqref{e1.4.30} implies that
$H(x)\in\dom D$ for $\mu$-a.e.\ $x\in \BX$. We therefore deduce
from Theorem \ref{t1.4.1} that
\begin{align*}
g(x,y):=D_y H(x)=\sum^\infty_{n=1} n I_{n-1}(h_n(x,y,\cdot))
\end{align*}
$\BP$-a.s.\ and for $\mu^2$-a.e.\ $(x,y)\in \BX^2$.
Using assumption \eqref{e1.4.30} together with the isometry properties
\eqref{e1.orth}, we infer that
\begin{align*}
\sum^\infty_{n=1} nn! \|\tilde {h}_n\|^2_{n+1}\le\sum^\infty_{n=1} n n! \|h_n\|^2_{n+1} 
= \BE \iint (D_yH(x))^2\mu(\dint x)\mu(\dint y) <\infty,
\end{align*}
yielding that $H\in\dom\delta$.

Now we define $H^{(m)}\in\dom\delta$, $m\in\N$, by \eqref{e1.4.20}
and note that
\begin{align*}
\BE \delta(H^{(m)})^2=\sum^m_{n=0}\BE I_{n+1}(\tilde{h}_n)^2
=\sum^m_{n=0}(n+1)! \|\tilde{h}_n\|^2_{n+1}.
\end{align*}
Using the symmetry properties of the functions $h_n$ it is easy
to see that the latter sum equals
\begin{align}\label{e1.4.34}
\sum^m_{n=0} n!\int h^2_n d\mu^{n+1}
+\sum^m_{n=1} n n!\iint h_n(x,y,z)
h_n(y,x,z)\mu^{2}(\dint(x,y))\mu^{n-1}(\dint z).
\end{align}
On the other hand, we have from Theorem \ref{t1.4.1} that
$$
D_y H^{(m)}(x)=\sum^m_{n=1} n I_{n-1}(h_n(x,y,\cdot)),
$$
so that
\begin{align*}
\BE\int H^{(m)}(x)^2\mu(\dint x)
+\BE \iint D_yH^{(m)}(x)D_xH^{(m)}(y)\mu(\dint x)\mu(\dint y)
\end{align*}
coincides with \eqref{e1.4.34}. Hence
\begin{align}\label{e1.4.36}
\BE \delta(H^{(m)})^2= \BE\int H^{(m)}(x)^2\mu(\dint x)
+\BE \iint D_yH^{(m)}(x)D_xH^{(m)}(y)\mu(\dint x)\mu(\dint y).
\end{align}
These computations imply that $g_m(x,y):=D_y H^{(m)}(x)$ converges in
$L^2(\BP\otimes\mu^2)$ towards $g$. Similarly,
$g'_m(x,y):=D_xH^{(m)}(y)$ converges towards $g'(x,y):=D_x g(y)$.
Since we have seen in the proof of Theorem  \ref{t1.4.3} that
$H^{(m)}\to H$ in $L^2(\BP\otimes\mu)$ as $m\to\infty$,
we can now conclude that the right-hand side of
\eqref{e1.4.36} tends to the right-hand side of the asserted
identity \eqref{e1.4.31}.
On the other hand we know by \eqref{e1.4.24} that
$\BE\delta(H^{(m)})^2\to \BE\delta(H)^2$
as $m\to\infty$. This concludes the proof.\qed

\bigskip

To explain the connection of \eqref{e1.4.30} with classical
stochastic analysis we assume for a moment that $\BX$ is equipped with
a transitive binary relation $<$ such that
$\{(x,y):x<y\}$ is a measurable subset of $\BX^2$
and such that $x<x$ fails for all $x\in \BX$.
We also assume
that $<$ totally orders the points of $\BX$ $\mu$-a.e., that is
\begin{align}\label{e1.4.27}
\mu([x])=0, \quad x\in \BX,
\end{align}
where $[x]:=\BX\setminus \{y\in \BX:\text{$y<x$ or $x<y$}\}$.
For any $\chi\in\bN_\sigma$ let $\chi_x$ denote the restriction
of $\chi$ to $\{y\in \BX:y<x\}$. Our final assumption on
$<$ is that $(\chi,y)\mapsto \chi_y$ is measurable.
A measurable function $h:\bN_\sigma\times \BX\to\R$ is called
{\em predictable} if
\begin{align}\label{e1.4.25}
h(\chi,x)=h(\chi_x,x),\quad (\chi,x)\in\bN_\sigma\times \BX.
\end{align}
A process $H\in L^0_\eta(\BP\otimes\mu)$ is predictable if it
has a predictable representative. In this case we have
$\BP\otimes\mu$-a.e.\ that $D_xH(y)=0$ for $y<x$ and
$D_yH(x)=0$ for $x<y$. In view of \eqref{e1.4.27} we obtain from
\eqref{e1.4.31} the classical It\^o isometry
\begin{align}\label{e1.4.35}
\BE \delta(H)^2= \BE\int H(x)^2\mu(\dint x).
\end{align}
In fact, a combinatorial argument shows that any predictable
$H\in L^2_\eta(\BP\otimes\mu)$ is in the domain of $\delta$.
We refer to \cite{LaPe11b} for more detail and references to
the literature.

We return to the general setting and derive a pathwise interpretation 
of the Kabanov-Skorohod integral. 
For $H\in L^1_\eta(\BP\otimes\mu)$ with representative $h$
we define
\begin{align}\label{pathSk}
\delta'(H):=\int h(\eta-\delta_x,x)\eta(\dint x)
-\int h(\eta,x)\mu(\dint x).
\end{align}
The Mecke equation \eqref{Mecke} implies that this definition
does $\BP$-a.s.\ not depend on the choice of the representative.
The next result (see \cite{LaPe11}) shows that the Kabanov-Skorohod integral 
and the operator $\delta'$ coincides on the intersection of their domains.
In the case of a diffuse intensity measure $\mu$ (and  requiring
some topological assumptions on $(\BX,\cX)$)
the result is implicit in \cite{Pic96b}.

\begin{theorem}\label{t1.4.4}
  Let $H\in L^1_\eta(\BP\otimes\mu)\cap \dom \delta$. Then
$\delta(H)=\delta'(H)$ $\BP$-a.s.
\end{theorem}
{\em Proof:} Let $H$ have representative $h$. The Mecke equation \eqref{Mecke} shows that 
$\BE\int |h(\eta - \delta_x,x)|\eta(\dint x)<\infty$
as well as
\begin{align}\label{adj3}
\BE\int D_xf(\eta)h(\eta,x)\mu(\dint x)=\BE f(\eta)\delta'(H),
\end{align}
whenever $f:\bN_\sigma\rightarrow \R$ is measurable and bounded.
Therefore we obtain from \eqref{partialint}
that $\BE F\delta'(H)=\BE F\delta(H)$ provided
that $F:=f(\eta)\in \dom D$. By Lemma \ref{l1.2.3}  the space of such
bounded random variables is dense in $L^2_\eta(\BP)$, so  we
may conclude that the assertion holds.\qed

\bigskip

Finally in this section we discuss the
{\em Ornstein-Uhlenbeck generator} $L$ whose domain is given by
all $F\in L^2_\eta$ satisfying
$$
\sum^\infty_{n=1}n^2 n! \|f_n\|^2_n<\infty.
$$
In this case one defines
\begin{align*}
LF:=- \sum_{n=1}^\infty n I_n(f_n).
\end{align*}
The (pseudo) {\em inverse} $L^{-1}$ of $L$ is given by
\begin{align}\label{e1.4.18}
L^{-1}F:=-\sum_{n=1}^\infty \frac{1}{n} I_n(f_n).
\end{align}
The random variable $L^{-1}F$ is well-defined for any  $F\in L^2_\eta$. 
Moreover, \eqref{e1.4.1} implies that $L^{-1}F\in\dom L$.
The identity $LL^{-1}F=F$, however, holds only if $\BE F=0$.

The three Malliavin operators $D,\delta$ and $L$ are connected by 
a simple formula:

\begin{proposition}\label{p1.4.6}
Let $F\in\dom L$. Then $F\in\dom D$, $DF\in \dom\delta$
and $\delta(DF)=-LF$.
\end{proposition} 
{\em Proof:} The relationship $F\in\dom D$ is a direct consequence
of \eqref{e1.4.1}. Let $H:=DF$. By Theorem \ref{t1.4.1} 
we can apply \eqref{e1.4.11} with $h_n:=(n+1)f_{n+1}$. We have
\begin{align*}
\sum^\infty_{n=0}(n+1)!\|h_n\|^2_{n+1}
=\sum^\infty_{n=0}(n+1)!(n+1)^2\|f_{n+1}\|^2_{n+1}.
\end{align*}
showing that $H\in\dom\delta$.
Moreover, since $I_{n+1}(\tilde{h}_n)=I_{n+1}(h_n)$ it follows that
\begin{align*}
\delta(DF)=\sum^\infty_{n=0}I_{n+1}(h_n)=\sum^\infty_{n=0}(n+1)I_{n+1}(f_{n+1})=-LF,
\end{align*}
finishing the proof.\qed

\bigskip

The following pathwise representation shows that
the Ornstein-Uhlenbeck generator can be interpreted as the
generator of a free {\em birth and death process} on $\BX$.

\begin{proposition}\label{p1.4.8}
Let $F\in\dom L$ with representative $f$ and assume
that $DF\in L^1_\eta(\BP\otimes\mu)$. Then
\begin{align}\label{e1.4.12}
LF=\int (f(\eta-\delta_x)-f(\eta))\eta(\dint x)
+\int (f(\eta+\delta_x)-f(\eta))\mu(\dint x).
\end{align}
\end{proposition} 
{\em Proof:} We use Proposition \ref{p1.4.6}.
Since $DF\in L^1_\eta(\BP\otimes\mu)$ we can apply
Theorem \ref{t1.4.4} and the result follows by a straigthforward
calculation.\qed

\section{Products of Wiener-It\^o integrals}\label{secproduct}

In this section we generalize Proposition \ref{pproduct} to products
of the form $I_p(f)I_q(g)$,
where  $f\in L^2(\mu^p)$ and $g\in L^2(\mu^q)$ for $p,q\in\N$.
To simplify the notation we assume that $f$ and $g$ are symmetric.
In this case we define for any $r\in\{0,\dots,p\wedge q\}$
(where $p\wedge q:=\min\{p,q\}$)
and $l\in[r]$ the {\em contraction} $f\ast^l_r g:\BX^{p+q-r-l}\rightarrow \R$
by
\begin{align}\label{e1.3.41}
f\ast^l_r g&(x_1,\dots,x_{p+q-r-l})\\ \notag
&:=\int f(y_1,\dots,y_l,x_1,\dots,x_{p-l})\times\\ \notag
&\qquad \times g(y_1,\dots,y_l,x_1,\dots,x_{r-l},x_{p-l+1},\dots,x_{p+q-r-l})
\mu^l(\dint (y_1,\dots,y_l)),
\end{align}
whenever these integrals are well-defined. 
In particular $f\ast^0_0 g=f\otimes g$.
In the case $q=1$ we have $f\ast^0_1 g=f\otimes^0_1 g$ and
$f\ast^1_1 g=f\otimes^1_1 g$; see \eqref{e1.2.81} and \eqref{e1.2.82}.

Under stronger integrability assumptions (and for diffuse intensity
measure), the following result 
has been proved in \cite{Surg84}. Our proof is quite different
and also independent of the proof of Proposition \ref{pproduct}.

\begin{proposition} \label{pproductgeneral}
Let $f\in L^2_s(\mu^p)$ and $f\in L^2_s(\mu^q)$ and assume that
$f\ast^l_r g\in L^2(\mu^{p+q-r-l})$ for
all $r\in\{0,\dots,p\wedge q\}$ and $l\in\{0,\dots,r-1\}$. 
Then
\begin{align}\label{e1.2.13}
I_p(f)I_q(g)=\sum^{p\wedge q}_{r=0} r!\binom{p}{r}\binom{q}{r}
\sum^r_{l=0}\binom{r}{l}I_{p+q-r-l}(f\ast^l_r g),\quad \BP\text{-a.s.}
\end{align}
\end{proposition}
{\em Proof:} We first note that the Cauchy-Schwarz inequality
implies $f\ast^r_r g\in L^2(\mu^{p+q-2r})$ for all
$r\in\{0,\dots,p\wedge q\}$.

We prove \eqref{e1.2.13} by induction on $p+q$. 
For $p\wedge q=0$ the assertion is trivial.
For the induction step we assume that $p\wedge q\ge 1$.
If $F,G\in L^0_\eta$, then an easy calculation (using
representatives) shows that
\begin{align}\label{prodrule}
D_x(FG)=(D_xF)G+F(D_xG)+(D_xF)(D_xG)
\end{align}
holds $\BP$-a.s.\ and for $\mu$-a.e.\ $x\in\BX$.
Using this together with Theorem \ref{t1.4.1}
we obtain that
\begin{align*}
D_x(I_p(f)I_q(g))=pI_{p-1}(f_x)I_q(g)
+qI_p(f)I_{q-1}(g_x)+pqI_{p-1}(f_x)I_{q-1}(g_x),
\end{align*}
where $f_x:=f(x,\cdot)$ and $g_x:=g(x,\cdot)$. 
We aim at applying the induction hypothesis to each of the
summands on the above right-hand side.
To do so, we note that
\begin{align*}
({f_x\ast^l_r g)}(x_1,\dots,x_{p-1+q-r-l})= 
f\ast^l_{r} g(x_1,\dots,x_{p-1-l},x,x_{p-1-l+1}\dots,x_{p-1+q-r-l})
\end{align*}
for all $r\in\{0,\dots,(p-1)\wedge q\}$ and $l\in\{0,\dots,r\}$
and
\begin{align*}
(f_x\ast^l_r g_x)(x_1,\dots,x_{p-1+q-1-r-l})= 
f\ast^l_{r+1} g(x,x_1,\dots,x_{p-1+q-1-r-l})
\end{align*}
for all $r\in\{0,\dots,(p-1)\wedge (q-1)\}$ and $l\in\{0,\dots,r\}$.
Therefore the pairs $(f_x,g)$,
$(f,g_x)$ and $(f_x,g_x)$ satisfy for
$\mu$-a.e.\ $x\in\BX$ the assumptions of the proposition.
The induction hypothesis implies that
\begin{align*}
D_x(I_p(f)I_q(g))=&\sum^{{(p-1)}\wedge q}_{r=0} r!p\binom{p-1}{r}\binom{q}{r}
\sum^r_{l=0}\binom{r}{l}I_{p+q-1-r-l}(f_x\ast^l_r g)\\
&+\sum^{{p}\wedge (q-1)}_{r=0} r!q\binom{p}{r}\binom{q-1}{r}
\sum^r_{l=0}\binom{r}{l}I_{p+q-1-r-l}(f\ast^l_r g_x)\\
&+\sum^{{(p-1)}\wedge (q-1)}_{r=0} r!pq\binom{p-1}{r}\binom{q-1}{r}
\sum^r_{l=0}\binom{r}{l}I_{p+q-2-r-l}(f_x\ast^l_r g_x).
\end{align*}
A straighforward calculation
(left to the reader) implies that the above right-hand side equals
\begin{align*}
\sum^{p\wedge q}_{r=0} r!\binom{p}{r}\binom{q}{r}
\sum^r_{l=0}\binom{r}{l}(p+q-r-l)I_{p+q-r-l-1}((\widetilde{f\ast^l_r g})_x),
\end{align*}
where the summand for $p+q-r-l=0$ has to be interpreted as $0$.
It follows that
\begin{align*}
D_x(I_p(f)I_q(g))=D_xG,\quad \BP\text{-a.s.},\,\mu\text{-a.e.\ $x\in\BX$},
\end{align*}
where $G$ denotes the right-hand side of \eqref{e1.2.13}.
On the other hand, the isometry properties \eqref{e1.orth} 
show that $\BE I_p(f)I_q(g)=\BE G$.
Since $I_p(f)I_q(g)\in L^1_\eta(\BP)$ we can use
the {\em Poincar\'e inequality} of Corollary \ref{c1.6.35}.
to conclude that
\begin{align*} 
\BE (I_p(f)I_q(g)-G)^2=0.
\end{align*}
This finishes the induction and the result is proved.
\qed

\bigskip

In the case $q=1$ equation \eqref{e1.2.13} says that
\begin{align}\label{e1.2.17}
I_p(f)I_1(g)=I_{p+1}(f\otimes g)+pI_p(f\otimes^0_1 g)+pI_{p-1}(f\otimes^1_1 g).
\end{align} 
This coincides with \eqref{e1.2.7} since
$p^{-1}\sum^p_{j=1}f\otimes^0_j g$ is the symmetrization
of the function $f\otimes^0_1 g$, while
$f\otimes^1_j g$ does not depend on $j$.
If $\{f\ne 0\}\subset B^p$ and $\{g\ne 0\}\subset B^q$ for some
$B\in\cX_0$ (as in Lemma \ref{l1.3.1}), 
then \eqref{e1.2.13} can be established by a direct
computation, just as in the proof of Proposition \ref{pproduct}.
The argument is similar to the proof of Theorem 3.1 in \cite{LaPeSchulThae14}.
The required integrability follows from
the Cauchy-Schwarz inequality; see \cite[Remark 3.1]{LaPeSchulThae14}.
In the case $q\ge 2$ we do not see, however, how to get from this special 
to the general case via approximation.

Equation \eqref{e1.2.13} can be further generalized so as 
to cover the case of a finite product of Wiener-It\^o integrals.
We again refer the reader to \cite{Surg84} as well as 
to \cite{TaPecc11,LaPeSchulThae14}.

\section{Mehler's formula}

In this section we aim at deriving a pathwise representation
of the inverse \eqref{e1.4.18} of the Ornstein-Uhlenbeck generator.
To give the idea we define for $F\in L^2_\eta$ with representation \eqref{chaos2}
\begin{align}\label{e1.5.1}
T_sF:=\BE F+\sum_{n=1}^\infty e^{-ns} I_n(f_n),\quad s\ge 0.
\end{align}
The family $\{T_s:s\ge 0\}$ is the {\em Ornstein-Uhlenbeck semigroup},
see e.g.\ \cite{Priv09} and also \cite{Nualart06} for the Gaussian case.  
If $F\in\dom L$ then it is easy to see that
$$
\lim_{s\to 0}\frac{T_sF-F}{s}=L
$$
in $L^2(\BP)$, see \cite[Proposition 1.4.2]{Nualart06} for the Gaussian case.
Hence $L$ can indeed be interpreted as the
generator of the semigroup. But in the theory of Markov processes it is well-known
(see e.g.\ the resolvent identities in \cite[Theorem 19.4]{Kallenberg}) 
that
\begin{align}\label{e1.5.2}
L^{-1}F=-\int^\infty_0 T_sF  ds,
\end{align}
at least under certain assumptions. What we therefore need is a pathwise
representation of the operators $T_s$. Our  guiding star is
the birth and death representation in Proposition \ref{p1.4.8}.

For $F\in L^1_\eta$ with representative $f$ we define,
\begin{align}\label{e1.5.3}
P_s F:=\int \BE[f(\eta^{(s)}+\chi)\mid\eta]\Pi_{(1-s)\mu}(\dint\chi),\quad s\in[0,1],
\end{align}
where $\eta^{(s)}$ is a {\em $s$-thinning} of $\eta$
 and where
$\Pi_{\mu'}$ denotes the distribution of a Poisson process with intensity
measure $\mu'$. The thinning  $\eta^{(s)}$ can be defined by
removing the points in \eqref{e1.1.6} independently of each other with probability $1-s$;
see \cite[p.\ 226]{Kallenberg}. 
Since
\begin{align}\label{e1.5.9}
\Pi_\mu=\BE\bigg[\int \I\{\eta^{(s)}+\chi\in\cdot\}\Pi_{(1-s)\mu}(\dint\chi)\bigg],
\end{align}
this definition does almost surely not depend
on the representative of $F$.
Equation \eqref{e1.5.9} implies in particular that
\begin{align}\label{e1.5.17}
\BE[P_sF]=\BE[F],\quad F\in L^1_\eta,
\end{align}
while Jensen's inequality implies for any $p\ge 1$ the contractivity property
\begin{align}\label{e1.5.18}
  \BE[(P_sF)^p]\le \BE[|F|^p],\quad  s\in[0,1],\, F\in L^2_\eta.
\end{align}

We prepare the main result of this section with
the folowing crucial lemma from \cite{LaPeSchu14}.

\begin{lemma}\label{l1.5.1} Let $F\in L^2_\eta$.
Then, for all $n\in\N$ and $s\in[0,1]$,
\begin{align}\label{e1.5.31} 
  D^n_{x_1,\dots,x_n}(P_s F)=s^nP_s D^n_{x_1,\dots,x_n}F, \quad
\mu^n\text{-a.e. } (x_1,\ldots,x_n)\in \BX^n,\;\BP\text{-a.s.}
\end{align}
In particular
\begin{align}\label{e1.5.32}
\BE [D^n_{x_1,\dots,x_n} P_s F] =s^n \BE [D^n_{x_1,\dots,x_n} F], \quad
\mu^n\text{-a.e.} (x_1,\dots,x_n)\in \BX^n.
\end{align}
\end{lemma}
{\em Proof:} 
To begin with, we assume that the representative of $F$
is given by $f(\chi)=e^{-\chi(v)}$ for some $v:\BX\rightarrow[0,\infty)$
such that $\mu(\{v>0\})<\infty$.
By the definition of a $s$-thinning,
\begin{equation}\label{e:ee}
\BE\big[e^{-\eta^{(s)}(v)}\mid \eta\big]=
\exp\bigg[\int \log\big((1-s)+s e^{-v(y)}\big)\eta(\dint y)\bigg],
\end{equation}
and it follows from Lemma 12.2 in \cite{Kallenberg} that
$$
\int \exp(- \chi(v))\Pi_{(1-s)\mu}(\dint\chi)
=\exp\bigg[-(1-s)\int (1-e^{-v})\dint\mu \bigg].
$$
Hence, the definition \eqref{e1.5.3} of the operator $P_s$ implies that
the following function $f_s$ is a representative of $P_sF$:
\begin{align*}
f_s(\chi):=\exp\bigg[-(1-s)\int \big(1-e^{-v}\big)\dint\mu\bigg]
\exp\bigg[\int \log\big((1-s)+s e^{-v(y)}\big)\chi(\dint y)\bigg].
\end{align*}
Therefore we obtain for any $x\in \BX$, that
\begin{align*}
D_xP_sF=f_s(\eta+\delta_x)-f_s(\eta)=s\big(e^{-v(x)}-1\big)f_s(\eta)
=s\big(e^{-v(x)}-1\big)P_sF.
\end{align*}
This identity can be iterated to yield for all $n\in\N$ and all
$(x_1,\ldots,x_n)\in\BX^n$ that
\begin{align*}
D^n_{x_1,\ldots,x_n}P_sF=s^n\prod^n_{i=1}\big(e^{-v(x_i)}-1\big)P_sF.
\end{align*}
On the other hand we have $\BP$-a.s.\ that
\begin{align*}
  P_sD^n_{x_1,\ldots,x_n}F=P_s\prod^n_{i=1}\big(e^{-v(x_i)}-1\big)F
= \prod^n_{i=1}\big(e^{-v(x_i)}-1\big)P_sF,
\end{align*}
so that \eqref{e1.5.31} 
holds for Poisson functionals of the given form.

By linearity, \eqref{e1.5.31} extends
to all $F$ with a representative in the set $\mathbf{G}$
of all linear combinations of functions $f$ as above.
There are $f_k\in \mathbf{G}$, $k\in\N$, satisfying
$F_k:=f_k(\eta)\to F=f(\eta)$ in $L^2(\BP)$ as $k\to\infty$,
where $f$ is a representative of $F$ (see \cite[Lemma 2.1]{LaPe11}).
Therefore
we obtain from the contractivity property \eqref{e1.5.18} that
\begin{align*}
\BE[ (P_sF_k-P_sF)^2]=\BE[(P_s(F_k-F))^2]\le \BE[(F_k-F)^2]\to 0,
\end{align*}
as $k\to\infty$. Taking $B\in\cX$ with $\mu(B)<\infty$,
it therefore follows from \cite[Lemma 2.3]{LaPe11} that
\begin{align*}
\BE\int_{B^n} |D^n_{x_1,\ldots,x_n}P_sF_k-D^n_{x_1,\ldots,x_n} P_sF|
\mu(\dint (x_1,\ldots,x_n))\to 0,
\end{align*}
as $k\to\infty$. On the other hand we obtain from the Fock 
space representation \eqref{e1.2.3} that $\BE |D^n_{x_1,\ldots,x_n}F|<\infty$
for $\mu^n$-a.e.\ $(x_1,\ldots,x_n)\in\BX^n$, so that
linearity of $P_s$ and \eqref{e1.5.18} imply
\begin{align*}
  \BE\int_{B^n} |P_sD^n_{x_1,\ldots,x_n}&F_k-P_sD^n_{x_1,\ldots,x_n}F|
\mu(\dint (x_1,\ldots,x_n))\\
&\le\int_{B^n} \BE |D^n_{x_1,\ldots,x_n}(F_k-F)|\mu(\dint (x_1,\ldots,x_n)).
\end{align*}
Again, this latter integral tends to $0$ as $k\to\infty$.
Since \eqref{e1.5.31} 
holds for any $F_k$ we obtain that
\eqref{e1.5.31} 
holds $\BP\otimes(\mu_B)^n$-a.e., and hence
also $\BP\otimes\mu^n$-a.e.

Taking the expectation in 
\eqref{e1.5.31} and using \eqref{e1.5.17} proves \eqref{e1.5.32}.
\qed

\bigskip

The following theorem from \cite{LaPeSchu14} 
achieves the desired pathwise representation of the
inverse Ornstein-Uhlenbeck operator.

\begin{theorem}\label{t1.5.2} Let $F\in L^2_\eta$. If $\BE F=0$
then we have $\BP$-a.s.\ that
\begin{align}\label{integratedMehler}
L^{-1}F=-\int^1_0 s^{-1}P_sF \dint s.
\end{align}
\end{theorem}
{\em Proof:}
Assume that $F$ is given as in \eqref{chaos2}.
Applying \eqref{chaos2} to $P_sF$ and using \eqref{e1.5.32} yields
\begin{equation}\label{e1.5.12}
P_sF=\BE F+\sum_{n=1}^\infty s^n I_n(f_n), \quad\BP\text{-a.s.},\, s\in[0,1].
\end{equation}
Furthermore,
$$
-\sum^m_{n=1}\frac{1}{n} I_n(f_n)=-\int^1_0s^{-1}\sum^m_{n=1}s^n I_n(f_n)\dint s,\quad
m\ge 1.
$$
Assume now that $\BE F=0$.
In view of \eqref{e1.4.18} we need to show that the above right-hand
side converges in $L^2(\BP)$, as $m\to\infty$, to the right-hand of side
of \eqref{integratedMehler}. Taking into account \eqref{e1.5.12} we hence have to show that
\begin{align*}
R_m:=\int^1_0s^{-1}\bigg(P_sF-\sum^m_{n=1}s^n I_n(f_n)\bigg) \dint s
=\int^1_0s^{-1}\bigg(\sum^\infty_{n=m+1}s^n I_n(f_n)\bigg) \dint s
\end{align*}
converges in $L^2(\BP)$ to zero. Using that
$\BE I_n(f_n)I_m(f_m)=\I\{m=n\}n!\|f_n\|^2_n$ we obtain
\begin{align*}
  \BE R^2_m\le\int^1_0s^{-2}\BE\bigg(\sum^\infty_{n=m+1}s^n I_n(f_n)\bigg)^2 \dint s
  =\sum^\infty_{n=m+1}n!\|f_n\|^2_n\int^1_0s^{2n-2} \dint s
\end{align*}
which tends to zero as $m\to\infty$.\qed

\bigskip

Equation \eqref{e1.5.12} implies {\em Mehler's formula}
\begin{align}\label{e1.5.16}
P_{e^{-s}}F=\BE F+\sum_{n=1}^\infty e^{-ns} I_n(f_n), \quad\BP\text{-a.s.},\,s\ge 0,
\end{align}
which was proved  in \cite{Priv09} for the special case of a finite Poisson process 
with a diffuse intensity measure.
Originally this formula was first established in a Gaussian
setting, see e.g.\ \cite{Nualart06}.
The family $\{P_{e^{-s}}:s\ge 0\}$ of operators describes
a special example of {\em Glauber dynamics}. 
Using \eqref{e1.5.16} in \eqref{integratedMehler} gives the identity \eqref{e1.5.1}.

\section{Covariance identities}

The fundamental Fock space isometry \eqref{e1.2.3} can be rewritten
in several other disguises. We give here two examples, 
starting with a covariance identity from 
\cite{HoP02} involving the operators $P_s$.

\begin{theorem}\label{t1.6.1} For any $F,G\in \dom D$,
\begin{align}\label{e1.6.2}
\BE[FG]=\BE[F]\BE[G]
+\BE\iint^1_0 (D_xF)(P_tD_xG) \dint t  \mu(\dint x).
\end{align}
\end{theorem}
{\em Proof:}
The Cauchy-Schwarz inequality
and the contractivity property \eqref{e1.5.18} imply that
\begin{align*}
\bigg(\BE\iint^1_0 |D_xF||P_sD_xG|\dint s \mu(\dint x)\bigg)^2
\le \BE\int (D_xF)^2\mu(\dint x) \BE\int (D_xG)^2\mu(\dint x)
\end{align*}
which is finite due to Theorem \ref{t1.4.1}. 
Therefore we can use Fubini's theorem
and \eqref{e1.5.31} to obtain that the right-hand side of \eqref{e1.6.2} equals
\begin{align}\label{e1.6.4}
  \BE[F]\BE[G]+\iint^1_0 s^{-1}\BE[(D_xF)(D_xP_sG)]\dint s\mu(\dint x).
\end{align}
For $s\in[0,1]$ and $\mu$-a.e.\ $x\in \BX$ we can apply the Fock space
isometry Theorem \ref{t1.2.1} to $D_xF$ and $D_xP_sG$.
Taking into account Lemma \ref{l1.5.1}, \eqref{e1.5.17} and
applying Fubini again (to be justified
below) yields that the second summand in \eqref{e1.6.4} equals
\begin{align*}
&\iint^1_0 s^{-1}\BE[D_xF]\BE[D_xP_sG]\dint s\mu(\dint x)\\
&+\sum^\infty_{n=1}\frac{1}{n!}\iiint^1_0
s^{-1}\BE[D^{n+1}_{x_1,\dots,x_n,x}F]\BE[D^{n+1}_{x_1,\dots,x_n,x}P_sG] \dint s\\
&\qquad\qquad\qquad\mu^n(\dint(x_1,\dots,x_n))\mu(\dint x)\\
=&\int \BE[D_xF]\BE[D_xG]\mu(\dint x)\\
&+\sum^\infty_{n=1}\frac{1}{n!}\iiint^1_0
s^{n}\BE[D^{n+1}_{x_1,\dots,x_n,x}F]\BE[D^{n+1}_{x_1,\dots,x_n,x}G] \dint s
\mu^n(\dint(x_1,\dots,x_n))\mu(\dint x)\\
=& \sum^\infty_{m=1}\frac{1}{m!}\int
\BE[D^m_{x_1,\dots,x_m}F]\BE[D^m_{x_1,\dots,x_m}G]\mu^m(\dint(x_1,\dots,x_m)).
\end{align*}
Inserting this into \eqref{e1.6.4} and applying Theorem \ref{t1.2.1}
yields the asserted formula \eqref{e1.6.2}. The use of Fubini's theorem
is justified by Theorem \ref{t1.2.1} for $f=g$  and the Cauchy-Schwarz inequality.
\qed

\bigskip

The integrability assumptions of Theorem \ref{t1.6.1} 
can be reduced to mere square integrability when using a symmetric formulation.
Under the assumptions of Theorem \ref{t1.6.1} the following result
was proved in \cite{HoP02}. An even more general version
is \cite[Theorem 1.5]{LaPe11}.

\begin{theorem}\label{t1.6.2} For any $F\in L^2_\eta$,
\begin{align}\label{e1.5.6}
\BE \iint^1_0 (\BE[D_xF\mid\eta^{(t)}])^2\dint t\mu(\dint x)<\infty,
\end{align}
and for any $F,G\in L^2_\eta$,
\begin{align}\label{e1.5.7}
\BE[FG]=\BE[F]\BE[G]+
\BE\iint^1_0 \BE[D_xF\mid\eta^{(t)}]\BE[D_xG\mid\eta^{(t)}]\dint t\mu(\dint x).
\end{align}
\end{theorem}
{\em Proof:} It is well-known (and not hard to prove)
that $\eta^{(t)}$ and $\eta-\eta^{(t)}$ 
are independent Poisson processes with
intensity measures $t\mu$ and $(1-t)\mu$, respectively.
Therefore we have for $F\in L^2_\eta$ with representative $f$ that
\begin{align}\label{e1.5.4}
\BE[D_xF|\eta_t]=\int D_x f(\eta^{(t)}+\chi)\Pi_{(1-t)\mu}(\dint\chi)
\end{align}
holds almost surely.
It is easy to see that the right-hand side
of \eqref{e1.5.4} is a measurable function of (the suppressed)
$\omega\in\Omega$, $x\in \BX$, and $t\in[0,1]$. 

Now we take $F,G\in L^2_\eta$ with representatives $f$ and $g$.
Let us first assume that $DF,DG\in L^2(\BP\otimes\mu)$.
Then \eqref{e1.5.6} follows from the (conditional) Jensen inequality
while \eqref{e1.5.4} implies for all $t\in[0,1]$ and $x\in \BX$, that
\begin{align*}
\BE[(D_xF) (P_t D_xG)]
&=\BE\bigg[D_xF \int D_xg(\eta^{(t)}+\mu)\Pi_{(1-t)\mu}(\dint\mu)\bigg]\\
&=\BE[\BE[D_xF\,\BE[D_xG\mid\eta^{(t)}]]=\BE[\BE[D_xF\mid\eta^{(t)}]\BE[D_xG\mid\eta^{(t)}]].
\end{align*}
Therefore \eqref{e1.5.7} is just another version of \eqref{e1.6.2}.

In this second step of the proof we consider general $F,G\in L^2_\eta$.
Let $F_k\in L^2_\eta$, $k\in\N$, be a sequence such that
$DF_k\in L^2(\BP\otimes\mu)$ and
$\BE (F-F_k)^2\to 0$ as $k\to\infty$.
We have just proved that
\begin{align*}
\Var [F_k-F^l]=
\BE\int(\BE[D_xF_k\mid\eta^{(t)}]-\BE[D_xF^l\mid\eta^{(t)}])^2\mu^*(\dint(x,t)), 
\quad k,l\in\N,
\end{align*}
where $\mu^*$ is the product of $\mu$ and Lebesgue measure on $[0,1]$.
Since $L^2(\BP\otimes\mu^*)$ is complete,
there is an $h\in L^2(\BP\otimes\mu^*)$ satisfying  
\begin{align}\label{345}
\lim_{k\to\infty}\BE\int(h(x,t)-\BE[D_xF_k\mid\eta^{(t)}])^2\mu^*(\dint(x,t))=0.
\end{align}
On the other hand it follows from Lemma \ref{lemsubs} that
for any $C\in\mathcal{X}_0$
\begin{align*}
\int_{C\times[0,1]} \BE\big|\BE[D_xF_k\mid\eta^{(t)}]&-\BE[D_xF\mid\eta^{(t)}]\big|
\mu^*(\dint(x,t))\\ 
&\le \int_{C\times[0,1]} \BE |D_xF_k-D_xF|\mu^*(\dint(x,t))\to 0
\end{align*}
as $k\to\infty$. Comparing this with \eqref{345} shows that
$h(\omega,x,t)=\BE[D_xF\mid\eta^{(t)}](\omega)$ for 
$\BP\otimes\mu^*$-a.e.\ $(\omega,x,t)\in\Omega\times C\times[0,1]$
and hence also for $\BP\otimes\mu^*$-a.e.\ 
$(\omega,x,t)\in\Omega\times \BX\times[0,1]$.
Therefore the fact that
$h\in L^2(\BP\otimes\mu^*)$ implies \eqref{e1.5.7}.
Now let $G_k$, $k\in\N$, be a sequence  approximating $G$.
Then equation  \eqref{e1.5.7} holds with $(F_k,G_k)$ instead of
$(F,G)$. But the second summand is just a scalar product
in $L^2(\BP\otimes\mu^*)$. Taking the limit as $k\to\infty$
and using the $L^2$-convergence proved above, yields
the general result. \qed

\bigskip

A quick consequence of the previous theorem  is the
{\em Poincar\'e inequality} for Poisson processes.
The following general version is taken from \cite{Wu00}.
A more direct approach can be based on the Fock space
representation in Theorem \ref{t1.2.1}, see \cite{LaPe11}. 

\begin{theorem}\label{t1.6.3} For any $F\in L^2_\eta$,
\begin{align}\label{Poincare}
\Var F\le \BE\int (D_xF)^2\mu(\dint x).
\end{align}
\end{theorem}
{\em Proof:} Take $F=G$ in \eqref{e1.5.7} and apply 
Jensen's inequality.\qed

\bigskip

The following extension of \eqref{Poincare} (taken from \cite{LaPeSchu14})
has been used in the proof of Proposition \ref{pproductgeneral}.

\begin{corollary}\label{c1.6.35}
For $F\in L^1_\eta$,
\begin{align}\label{eq:PoincareL1}
\BE F^2 \leq (\BE F)^2 + \BE \int (D_xF)^2\mu(\dint x).
\end{align}
\end{corollary}
{\em Proof:}
For $s>0$ we define
$$
F_s=\I\{F>s\}s+ \I\{-s\leq F\leq s\}F -\I\{F<-s\}s 
$$
By definition of $F_s$ we have $F_s\in L^2_\eta$ and $|D_xF_s|\leq |D_xF|$ 
for $\mu$-a.e.\ $x\in\BX$. Together with the Poincar\'e inequality
\eqref{Poincare} we obtain that
$$
\BE F_s^2 \le (\BE F_s)^2
+\BE\int (D_xF_s)^2 \mu(\dint x)\leq  (\BE F_s)^2+\BE\int (D_xF)^2 \mu(\dint x).
$$
By the monotone convergence theorem and 
the dominated convergence theorem, respectively, we have that 
$\BE F_s^2\to\BE F^2$ and $\BE F_s\to\BE F$ as $s\to\infty$. 
Hence letting $s\to\infty$ in the previous inequality yields the assertion.
\qed

\bigskip

As a second application of Theorem \ref{t1.6.2} we obtain
the {\em Harris-FKG inequality} for Poisson processes, derived
in \cite{Janson84}. Given $B\in \cX$,
a function $f\in \mathbf{F}(\bN_\sigma)$ is {\em increasing on} $B$
if $f(\chi+\delta_x)\ge f(\chi)$ for all $\chi\in \bN_\sigma$
and all $x\in B$. It is  
{\em decreasing on} $B$ if $(-f)$ is increasing on $B$.

\begin{theorem}\label{t1.6.4}
Suppose $B \in \cX$.
Let $f,g\in L^2(\BP_\eta)$ be increasing on $B$ and decreasing on
$\BX \setminus B$. Then
\begin{align}\label{FKG}
\BE [f(\eta)g(\eta)]\ge (\BE f(\eta))(\BE g(\eta)).
\end{align}
\end{theorem}

It was noticed in \cite{Wu00} that the correlation inequality \eqref{FKG}
(also referred to as {\em association})
is a direct consequence of a covariance identity. 

\bigskip

\noindent {\bf Acknowledgment:}
The proof of Proposition \ref{pproductgeneral} is joint work
with Matthias Schulte.

\end{document}